
\def\LaTeX{%
  \let\Begin\begin
  \let\End\end
  \let\salta\relax
  \let\finqui\relax
  \let\futuro\relax}

\def\UK{\def\our{our}\let\sz s}
\def\USA{\def\our{or}\let\sz z}

\UK



\LaTeX

\USA


\salta

\documentclass[twoside,12pt]{article}
\setlength{\textheight}{24cm}
\setlength{\textwidth}{16cm}
\setlength{\oddsidemargin}{2mm}
\setlength{\evensidemargin}{2mm}
\setlength{\topmargin}{-15mm}
\parskip2mm


\usepackage[usenames,dvipsnames]{color}
\usepackage{amsmath}
\usepackage{amsthm}
\usepackage{amssymb}
\usepackage[mathcal]{euscript}
\usepackage{enumitem}
%
%
\usepackage{cite}
%
%
%


\definecolor{viola}{rgb}{0.3,0,0.7}
\definecolor{ciclamino}{rgb}{0.5,0,0.5}
\definecolor{rosso}{rgb}{0.8,0,0}

\def\dafare #1{{\color{red}#1}}

\def\dafare #1{#1}


\bibliographystyle{plain}


%

\finqui

\def\Beq{\Begin{equation}}
\def\Eeq{\End{equation}}
\def\Bsist{\Begin{eqnarray}}
\def\Esist{\End{eqnarray}}

\def\Bthm{\Begin{theorem}}
\def\Ethm{\End{theorem}}
\def\Blem{\Begin{lemma}}
\def\Elem{\End{lemma}}
\def\Bprop{\Begin{proposition}}
\def\Eprop{\End{proposition}}

\def\Bdim{\Begin{proof}}
\def\Edim{\End{proof}}
\def\Bcenter{\Begin{center}}
\def\Ecenter{\End{center}}
\let\non\nonumber




\def\step #1 \par{\medskip\noindent{\bf #1.}\quad}


\def\Lip{Lip\-schitz}
\def\Holder{H\"older}
\def\Frechet{Fr\'echet}

\def\aand{\quad\hbox{and}\quad}

\def\wk{well-known}
\def\socal{so-called}
\def\lhs{left-hand side}
\def\rhs{right-hand side}
\def\sfw{straightforward}

\def\CH{Cahn-Hilliard}



\def\multibold #1{\def\arg{#1}%
  \ifx\arg\pto \let\next\relax
  \else
  \def\next{\expandafter
    \def\csname #1#1#1\endcsname{{\bf #1}}%
    \multibold}%
  \fi \next}

\def\pto{.}

\def\multical #1{\def\arg{#1}%
  \ifx\arg\pto \let\next\relax
  \else
  \def\next{\expandafter
    \def\csname cal#1\endcsname{{\cal #1}}%
    \multical}%
  \fi \next}


\def\multimathop #1 {\def\arg{#1}%
  \ifx\arg\pto \let\next\relax
  \else
  \def\next{\expandafter
    \def\csname #1\endcsname{\mathop{\rm #1}\nolimits}%
    \multimathop}%
  \fi \next}

\multibold
qwertyuiopasdfghjklzxcvbnmQWERTYUIOPASDFGHJKLZXCVBNM.

\multical
QWERTYUIOPASDFGHJKLZXCVBNM.

\multimathop
ad dist div dom meas sign supp .


\def\accorpa #1#2{\eqref{#1}-\eqref{#2}}
\def\Accorpa #1#2 #3 {\gdef #1{\eqref{#2}-\eqref{#3}}%
  \wlog{}\wlog{\string #1 -> #2 - #3}\wlog{}}


\def\graffe #1{\mathopen\{#1\mathclose\}}

\def\<#1>{\mathopen\langle #1\mathclose\rangle}
\def\norma #1{\mathopen \| #1\mathclose \|}

\def\intQt{\int_{Q_t}}
\def\intQtT{\int_{Q_t^T}}
\def\intQ{\int_Q}
\def\iO{\int_\Omega}

\def\inttt{\int_t^T}

\def\dt{\partial_t}
\def\dn{\partial_n}
\def\dtt{\partial_{tt}}

\def\checkmmode #1{\relax\ifmmode\hbox{#1}\else{#1}\fi}

\def\aeQ{\checkmmode{a.e.\ in~$Q$}}

\def\aaQ{\checkmmode{for a.a.~$(x,t)\in Q$}}

\def\aat{\checkmmode{for a.a.~$t\in(0,T)$}}


\def\erre{{\mathbb{R}}}

\def\enne{{\mathbb{N}}}




\def\genspazio #1#2#3#4#5{#1^{#2}(#5,#4;#3)}
\def\spazio #1#2#3{\genspazio {#1}{#2}{#3}T0}

\def\L {\spazio L}
\def\H {\spazio H}

\def\C #1#2{C^{#1}([0,T];#2)}


\def\Lx #1{L^{#1}(\Omega)}
\def\Hx #1{H^{#1}(\Omega)}

\def\Luno{\Lx 1}
\def\Ldue{\Lx 2}

\def\Huno{\Hx 1}
\def\Hdue{\Hx 2}



\let\theta\vartheta
\let\eps\varepsilon
\let\phi\varphi

\let\TeXchi\chi                         
\newbox\chibox
\setbox0 \hbox{\mathsurround0pt $\TeXchi$}
\setbox\chibox \hbox{\raise\dp0 \box 0 }
\def\chi{\copy\chibox}



\def\bQ{b_1}

\def\bQh{b_2}
\def\bOh{b_3}
\def\bz{b_0}

\def\phQ{\phi_Q}

\def\phO{\phi_\Omega}


\def\sO{\s_\Omega}

\def\sQ{\s_Q}

\def\Uad{\calU_{\ad}}
\def\uopt{\overline u}

\def\Vp{V^*}

\def\normaV #1{\norma{#1}_V}
\def\normaH #1{\norma{#1}_H}

\let\hat\widehat

\def\Pi{\hat\pi}


\def\cd{c_\delta}
\def\s{\sigma}  
\def\m{\mu}	    
\def\ph{\phi}	
\def\a{\alpha}	
\def\b{\beta}	
\def\d{\delta}  
\def\et{\eta}   
\def\bph{\overline\ph}  
\def\bm{\overline\m}    
\def\bs{\overline\s}    
\def\J{{\cal J}}

\def\I2 #1{\int_{Q_t}|{#1}|^2}
\def\IN2 #1{\int_{Q_t}|\nabla{#1}|^2}
\def\IO2 #1{\iO |{#1(t)}|^2}
\def\INO2 #1{\iO |\nabla{#1}(t)|^2}
\def\UR{{\cal U}_R}
\def\CP{\boldsymbol{(CP)}}

\def\qab {q_{\b}}
\def\pab {p_{\b}}
\def\rab {r_{\b}}


\Begin{document}

\title{
		Optimal treatment for a phase field system
		\\[0.3cm] 
		of Cahn-Hilliard type  
		modeling tumor growth 		
		\\[0.3cm] 
		by asymptotic scheme }
\author{}
\date{}
\maketitle

\Bcenter
\vskip-1cm
{\large\sc Andrea Signori$^{(1)}$}\\
{\normalsize e-mail: {\tt andrea.signori02@universitadipavia.it}}\\[.25cm]
$^{(1)}$
{\small Dipartimento di Matematica e Applicazioni, Universit\`a di Milano-Bicocca}\\
{\small via Cozzi 55, 20125 Milano, Italy}

\Ecenter
\Begin{abstract}\noindent
We consider a particular phase field system which physical context
is that of tumor growth dynamics. The model we deal with consists of a Cahn-Hilliard
equation governing the evolution of the phase variable which 
takes into account the tumor cells proliferation in the tissue
coupled with a reaction-diffusion equation for the nutrient.
This model has already been investigated
from the viewpoint of well-posedness, long-time behavior, and asymptotic analyses
as some parameters go to zero.
Starting from these results, we aim to face a related optimal control problem 
by employing suitable asymptotic schemes. 
In this direction, we assume some quite general growth conditions for 
the involved potential and a smallness
restriction for a parameter appearing in the system we are going to face.
We provide existence of optimal controls and a necessary condition
for optimality is addressed.

\End{abstract}

\vskip3mm

\pagestyle{myheadings}
\newcommand\testodispari{\sc Signori}
\newcommand\testopari{\sc {Asymptotic analyses of a control problem}}
\markboth{\testodispari}{\testopari}

\salta
\finqui

\newpage
\section{Introduction}
\label{INTRODUCTION}
\setcounter{equation}{0}
Over the last decades, there has been increasing attention by the mathematical
community towards biological models for tumor growth (see \cite{CL}). 
Among them, the ones introduced by exploiting phase field approaches
and continuum mixture theory cover an important role.
The key idea consists in
reading the physical evolution process like an interaction between two particular
fluids which are designed to model the tumor cells and the healthy ones.
In this regards, we especially point out two classes.
The first one gives rise to the \socal~Cahn-Hilliard-Darcy or 
Cahn-Hilliard-Brinkman systems
which describe the tumor and healthy cells as inertia-less 
fluids including effects generated by fluid flow development.
To this concerns, let us refer to \cite{WLFC,GLSS,DFRGM ,GARL_1,GARL_4,GAR}.
The second class neglects the velocity and consists of a Cahn-Hilliard equation
(see, e.g., \cite{Mir_CH} and the huge references therein) for the phase 
variable coupled with a reaction-diffusion equation for the nutrient.
Moreover, let us point out the papers \cite{GARL_2, GARL_3}, 
where transport mechanisms such as chemotaxis and active transport are also taken into account.
Further investigations and mathematical models
related to biology can be found, e.g., in \cite{FLRS,FRL}.

Here, we try to describe the main purpose of the work postponing
the technicalities and the investigation of 
the proper assumptions that will be specified in the forthcoming section.
Before moving on, let us mention that with $\Omega \subset \erre^3$ we denote the 
set where the evolution takes place and, for a given final time $T>0$, we define
the standard parabolic cylinder and its boundary by
\Bsist
	& \non
	Q_{t}:=\Omega \times (0,t), \quad \Sigma_{t}:=\partial\Omega \times (0,t)
	\quad \hbox{for every $t \in (0,T]$,}
	\\
	& Q:=Q_{T}, \quad \hbox{and} \quad \Sigma:=\Sigma_{T}.
	\label{QS}
\Esist
Hence, the model we are going to consider reads as follows
\Bsist
  & \a \dt \m_{\b} + \dt \ph_{\b} - \Delta \m_{\b} = P(\phi_{\b}) (\s_{\b} - \m_{\b})
  \quad \hbox{in $\, Q$}
  \label{EQabprima}
  \\
  & \m_{\b} = \beta \dt \ph_{\b} - \Delta \ph_{\b} + F'(\ph_{\b})
  \label{EQabseconda}
  \quad \hbox{in $\,Q$}
  \\
  & \dt \s_{\b} - \Delta \s_{\b} = - P(\ph_{\b}) (\s_{\b} - \m_{\b}) + u_\b
  \label{EQabterza}
  \quad \hbox{in $\,Q$}
  \\
  & \dn \m_{\b} = \dn \ph_{\b} =\dn \s_{\b} = 0
  \quad \hbox{on $\,\Sigma$}
  \label{BCEQab}
  \\
  & \m_{\b}(0)=\m_0,\, \ph_{\b}(0)=\ph_0,\, \s_{\b}(0)=\s_0
  \quad \hbox{in $\,\Omega,$}
  \label{ICEQab}
\Esist
\Accorpa\EQab EQabprima ICEQab
for some positive constants $\a$ and $\b$.
Let us emphasize that the notation $\ph_\b$,
instead of the simplest $\ph$, and the same goes
for the other variables, is motivated by the fact that in the following
we are going to let $\b\searrow 0$ and we will denote as
$\ph$ the limit of $\ph_\b$. So, with the subscript $\b$, we aim to stress the
fact that such a variable corresponds to the system with $\b>0$.

The above system is a simplified version of the diffuse interface model originally 
introduced by Hawkins-Daruud et al. in \cite{HDZO}, where
the velocity and chemotaxis contributions are neglected (see also \cite{HKNZ,WZZ,CLLW,HDPZO,OHP})
and it also includes some regularizing parameters.
It is worth spending some words explaining the physical meaning of the model.
The unknown $\ph_\b$ is an order parameter and it is devoted to keeping track of
the evolution of the tumor in the tissue. 
It is usually normalized between $-1$ and $+1$, where 
these extreme values represent the pure phases,
that is the tumor phase and the healthy cell phase, respectively. 
The second unknown $\m_\b$, as usual for \CH~equation, stands
for the chemical potential for $\ph_\b$.
Finally, the last unknown $\s_\b$ represents the nutrient-rich extra-cellular water volume fraction.
It takes values between $0$ and $1$ with the following property:
the closer to one, the richer of water the extra-cellular fraction is, while the closer to zero,
the poorer it is.
Furthermore, $u_\b$ is the \socal~control variable which will allow us to
interact in some sense with the above system. As usual in control theory,
the system in which the control variable appears is referred to as state system.
Likewise, we will refer to the solution $(\m_\b,\ph_\b,\s_\b)$ to as the associated state.
As far as $P$ and $F$ are concerned, they are nonlinearities. The former is a proliferation function,
while the latter is a double-well potential. Customary examples for $F$ are the classical 
regular potential, the logarithmic potential, and the double-obstacle one.
In this contribution, we focus the attention on the first one which is given by
\Bsist
  & F_{reg}(r) := \frac14(r^2-1)^2 
  = \frac 14 ((r^2 - 1)^+)^2 + \frac 14 ((1 - r^2)^+)^2 \, 
  \quad \hbox{for } r \in \erre.
  \label{F_reg}
\Esist
For different physically meaningful choices of the 
potentials we refer to \cite{Agosti}, and
to the references therein, where several numerical applications to tumor growth
can be found as well.

We point out that the above model has been quite well-understood owing to the previous works 
\cite{CGH, CGRS_ASY, CGRS_VAN} as far as well-posedness and
long-time behavior are concerned. 
Moreover, in \cite{FGR} the analyses 
of the same model without the relaxation terms
$\a\dt\m_\b$ and $\b\dt\ph_\b$, has been performed considering regular potentials
and allowing $P$ to possess polynomial growth.
Besides, as long-time behavior of solutions are concerned, we also mention  
\cite{Kur}, where the author extends the well-posedness
results proved in \cite{CGRS_ASY,CGRS_VAN}, as $\b\searrow 0$,
to the case of unbounded domains. Moreover, we refer to
\cite{MRS}, where the authors investigate the long-time behavior of the
non-relaxed version of system \EQab, i.e. the case $\a=\b=0$, 
showing the existence of the global attractor in proper phase space.
In view of such flourishing literature, a further natural aim is to 
investigate some corresponding optimal control problems in which the state system is 
governed by the evolution system \EQab.
In this direction, we refer to the recent work \cite{S}, where, making extensive use of
the terms $\a \dt \m_\b$ and $\b \dt \ph_\b$, an optimal control problem for the above
system is tackled in a general framework for the potential, allowing both the classical
and the logarithmic potential to be considered. 
Additionally, the same author proves in the subsequent work \cite{S_DQ}, via a 
proper asymptotic scheme
known in the literature as to deep quench limit, that it is also possible to generalize
the assumptions for the potentials in order to take into account also singular and 
nonregular potentials like the double-obstacle.
Furthermore, let us refer to \cite{CGRS_OPT}, where
a similar optimal control problem is considered for the state 
system \EQab~without these relaxation terms.
In addition, regarding some optimal control problem in which time is taken into account,
we address the recent \cite{CRW}, where also 
long-time behavior of solutions
has been investigated, and we also mention \cite{GARLR}, 
where an optimal time therapeutic treatment has been investigated.
Finally, we point out the contributions \cite{EK,EK_ADV} in which
the authors investigate an optimal control problem for different tumor growth 
model based on the Cahn-Hilliard-Brinkman equation, which was previously
investigated in \cite{EGAR}, pointing out also some sufficient
conditions for the optimality.

Here, we actually aim to study an intermediate optimal control problem
between \cite{CGRS_OPT} and \cite{S}. In fact,
we still consider the state system to be \EQab, but without considering the
relaxation term $\b \dt \ph_\b$. 
\dafare{
Let us emphasize the mathematical interest of this problem.
On the one hand, the non-relaxed version of system \EQab~was 
investigated in \cite{FGR} and the corresponding
control problem was then tackled in \cite{CGRS_OPT}. 
On the other hand, the relaxed model
was studied in \cite{CGH}. There, due to the stronger regularity, the authors can 
take into account very general potentials. Then, by \cite{CGRS_ASY,CGRS_VAN},
the asymptotic analysis as $\a$ and $\b$ go to zero was performed.
Moreover, the corresponding control problem was recently treated by the author in \cite{S}.
Hence, it will be interesting to understand how the optimal
control problems associated with these similar models behave.
Furthermore, we refer to \cite{S_a}, where the author, focus the attention on the case
$\a\searrow 0, \b>0$.}

As for the control problem, we are going to take into account the following
tracking-type cost functional 
\begin{align}
	&
	\J (\ph, \s, u)  : = 
	\frac \bQ 2 \norma{\ph - \phQ}_{L^2(Q)}^2
	+ \frac \bQh 2 \norma{\s - \sQ}_{L^2(Q)}^2
	+ \frac \bOh 2 \norma{\s(T)-\sO}_{L^2(\Omega)}^2
    + \frac \bz 2 \norma{u}_{L^2(Q)}^2,
    \label{costfunct}
\end{align}
and the control-box constraints
\Beq
    \label{Uad}
	\Uad := \graffe{u \in L^{\infty}(Q): u_* \leq u \leq u^*\ \aeQ},
\Eeq
where $u_*$ and $u^*$ are functions that fix the admissible set in which the control variable 
$u$ can be chosen.
Furthermore, $\bz, \bQ, \bQh,\bOh$ stand for nonnegative constants, not all zero,
while $\phQ, \sQ, \sO $ denote some target functions defined in $Q$ and 
$\Omega$, respectively.

Since our starting point is \cite{S}, we will refer several times to the results 
there proved. 
So, it is worth noting that the cost functional \eqref{costfunct}
is slightly less general with respect to the one there proposed. 
There, an additional term of the form $\frac k 2 \norma{\ph(T)-\phO}_{L^2(\Omega)}^2$
also appears, for a nonnegative constant $k$ and 
a prescribed function $\phO$
which models the final configuration of the tumor colony. 
From a control viewpoint, this contribution allows us to force 
the final configuration of the tumor to be as close as possible, 
in the sense of $\Lx2$-norm, to the fixed target $\phO$.
Here, we restrict the investigation to the case $k=0$.
This will be motivated by the analysis of the corresponding adjoint problem
that leads to assume such a compatibility condition. 
\dafare{
To motivate this limitation, let us recall that (see \cite[Syst.~(2.22)-(2.26)]{S})
the final conditions of the adjoint problem pointed out in \cite{S}
are the following (note that the constant $k$ is called $b_2$ in that paper)
\Beq
	\non
	p(T) - \b q(T) = k (\bph(T)-\phO), \ \ \a p(T)=0.
\Eeq
However, if $\a>0$ and $\b=0$, we formally deduce that 
\Beq
	\non
	\begin{cases}
	 p(T) = k (\bph(T)-\phO)
     \\  
     \a p(T) =0,
	\end{cases}
\Eeq
which yields $p(T)=0$ and therefore also that $0= k (\bph(T)-\phO)$,
which is not satisfied in general since $\phO$ is arbitrary.}
Hence, to not lead to a contradiction we assume $k=0$, 
so that the choice of the less general cost functional \eqref{costfunct}
is now justified.

Thus, the optimal control treated in \cite{S} consists of 
solving the problem:
\Bsist
	\non
	\CP_{\boldsymbol{\b}} && \hbox{Minimize $\J(\ph,\m,u)$}
	\hbox{ subject to the control contraints \eqref{Uad} and under the}
	\\ && \non
	\hbox{requirement that the variables $(\ph, \s)$ 
	yield a solution to \EQab.}
	\label{P_ab}
\Esist
There, the author confirmed the existence of, at least, one optimal control
and also provide some first-order optimality conditions reading as
variational inequalities. 

Moreover, let us recall that the asymptotic analysis for 
system \EQab~has already been investigated in \cite{CGH, CGRS_ASY, CGRS_VAN}, 
where the authors carefully point out some sufficient conditions to let 
$\a$ and $\b$ go to zero, both sequentially and separately. 
As a matter of fact, they proved that, as $\b\searrow 0$, which is 
the case we are going to consider, providing to require additional assumptions,
the unique solution to system \EQab~converges to some limit which yields a solution to the
following problem
\Bsist
  & \non {}_{\Vp} \< \dt (\a \m + \ph) ,  v>  _V
  + \iO \nabla \m \cdot \nabla v 
  = \iO P(\phi) (\s - \m)v
  \\ &  \quad \hbox{$\forall v \in V$, $a.e.$ in $(0,T)$}
  \label{EQaprima}
  \\
  & \m = - \Delta \ph + F'(\ph)
  \label{EQaseconda}
  \quad \hbox{in $\,Q$}
  \\
  & \dt \s - \Delta\s = - P(\ph) (\s - \m) + u
  \label{EQaterza}
   \quad \hbox{in $\,Q$}
  \\
  & \dn \m = \dn \ph =\dn \s = 0
  \quad \hbox{on $\,\Sigma$}
  \label{BCEQa}
  \\
  & (\a \m + \ph) (0)=\a \m_0 + \ph_0, \,\, \s(0)=\s_0
  \quad \hbox{in $\,\Omega.$}
  \label{ICEQa}
\Esist
\Accorpa\EQa EQaprima ICEQa
Furthermore, it was shown that, under a suitable smallness requirement on $\a$,
that will be precised later on,
the solution is indeed unique.
%

From a different perspective, one can take the above system as a starting point, 
trying to face the analysis of the corresponding optimal control problem.
Namely, one can try to solve the following minimization problem:
minimize $\J (\ph, \s, u)$ subject to the control constraints \eqref{Uad}
and under the requirement that the variables $(\ph, \s)$ are solutions to \EQa.
Actually, this is the optimal control problem we try to tackle by following a different 
strategy consisting of pass
to the limit, as $\b \searrow 0$, in the optimal control problem $(CP)_{{\b}}$.
This technique turns out to be particularly interesting since
we still will obtain similar results with respect to \cite{S}. 
At the same time, we will treat the optimal control problem avoiding the investigation
of the linearized system, which is usually not so difficult
and, mostly, we can avoid
the discussion on the \Frechet~differentiability of the associated control-to-state mapping,
which is usually more challenging.

On the other hand, the first-order necessary condition of
$(CP)$ cannot be directly obtained by letting $\b\searrow 0$ in the optimality 
condition for the corresponding optimal control problem with $\b>0$.
This would be the case if we ensure that
every optimal control for $(CP)$ can be recovered as limit of sequence of 
optimal controls for $(CP)_\b$, which is quite a strong requirement.
\dafare{
Unfortunately, we are unable to prove such a global result.
However, a partial one can be stated localizing the problem by following the idea
firstly introduced by Barbu in {\cite{BARBU}} (see also, e.g., 
\cite{S_DQ,CFGS,CFS,CGS_DQ},
where such a technique was applied).
The idea consists in locally perturbing the problem $(CP)_\b$ in order to 
find the desired approximation result.
For this purpose, the main ingredient is the \socal~adapted cost functional that
is defined as follows
\Beq
	\label{adapted_cost}
	\widetilde{\J} (\ph,\s,u):=
	\J (\ph,\s,u)
	+ \frac 12 \norma{u-\overline u}_{L^2(Q)}^2,
\Eeq
where $\overline{u}$ stands for an optimal control for $(CP)$.
Thus, we are naturally lead to solving the following minimization problem:
\Bsist
	\non
	(\widetilde{CP})_{\b} && \hbox{Minimize $\widetilde{\J}(\ph,\m,u)$}
	\hbox{ subject to the control constraints \eqref{Uad} and under the}
	\\ && \non
	\hbox{requirement that the variables $(\ph, \s)$ 
	yield a solution to \EQab.}
	\label{Ptilde}
\Esist
Next, instead of looking for approximating sequence of optima for $(CP)_\b$,
one take a sequence of controls which are optimal for $(\widetilde{CP})_{\b}$.}
This technique will allow us to properly
let $\b\searrow 0$ in the optimality condition
for $(\widetilde{CP})_{\b}$ to recover the variational 
inequality which characterizes the necessary condition for optimality of $(CP)$.

Summing up, the current contribution has the purpose of showing, 
by means of asymptotic approaches,
that the following control problem admits a solution:
\Bsist
	\non
	\CP && \hbox{Minimize $\J(\ph,\m,u)$}
	\hbox{ subject to the control constraints \eqref{Uad} and under the}
	\\ && \non
	\hbox{requirement that the variables $(\ph, \s)$ 
	yield a solution to \EQa.}
	\label{P_a}
\Esist
Moreover, we will also provide a necessary condition that an optimal control has
to satisfy in terms of a variational inequality.

Lastly, let us sketch the physical background of
the control problem we are dealing with. Roughly speaking, we
are looking for the best choice $u$ in such 
a way that, with the corresponding solution to \EQa, it minimizes the cost 
functional $\J$. 
Furthermore, the control $u$ appears in equation \eqref{EQaterza}, the one describing
the evolution of the nutrient. Therefore, from the model viewpoint, it can 
be read as a supply of a nutrient or a drug in medical treatment. 
Moreover, for some given a priori targets $\phQ,\sQ,\sO$,
minimizing the cost functional $\J$ corresponds to force the system to approach a prescribed
targets which should be taken as desirable configurations
for clinical reasons, e.g., for surgery.
In addition, the ratios among the constants $\bz,\bQ,\bQh,\bOh$
implicitly describe which targets hold the leading part in our application
and the last term of the cost functional represents the cost we have to pay to 
take $u$ into account. In fact, it should be read as the
rate of risks to afflict harm to the patient by following that strategy,
namely the side-effect that may occur if too many drugs are dispensed
to the patient. 

The plan for the rest of the paper is as follows.
In Section~2, we set the notation we are going to use and recollect
the obtained results. From Section~3 onward, we start with the proofs of the stated results.
Furthermore, Section~3 is devoted to investigating the well-posedness and 
the asymptotic behavior, as $\b\searrow 0$, of the state system. 
Lastly, in Section~4, we discuss the control problem $(CP)$ by invoking some 
asymptotic schemes. We check the existence of optimal control, 
study the well-posedness
of the adjoint system, and provide first-order necessary condition
for optimality.
\section{Assumptions and Main Results}
\label{RESULTS}
\setcounter{equation}{0}

In this section, we aim to set the notation and collect our results.
To begin with, we recall that $\Omega$ stands for the space
domain where the evolution
takes place and we assume it to be a bounded, connected, smooth, and open set of $\erre^3$,
with boundary indicated by $\Gamma$.
As the functional spaces are concerned, it turns out to be 
convenient to introduce the following
\Beq
	\non
	H:= \Ldue, \quad V:= \Huno, \quad W:=\graffe{v \in \Hdue : \dn v = 0 \hbox{ on } \Gamma},
\Eeq
where $\partial _n$ stands for the outward normal derivative.
Furthermore, to work with Banach spaces, we
endow them with their standard norms.
For a general Banach space $X$, we use $\norma{\cdot}_{X}$ to designate its norm,
$X^*$ for its topological dual, and 
${}_{X^*}\<\cdot,\cdot>_{X}$ for the duality product between $X^*$ and $X$.
Likewise, we use the symbol $\norma{\cdot}_{p}$ for the usual norm in $L^p(\Omega)$.
The above definitions, in turn, imply that $(V,H,\Vp)$ constitutes a Hilbert triplet,
that is, the following injections $V \subset H \equiv H^* \subset V^*$ 
are both continuous and dense and we have the standard identification
between the dual product of $V$ and the inner product of $H$.
Namely, it holds that 
\Beq
	\non
	{}_{\Vp}\< u,v >_V = \iO u v \quad \hbox{for every $u \in H$ and $v \in V$.}
\Eeq
As for the basic assumptions for the system \EQab~and for the cost functional 
\eqref{costfunct}, we postulate that
\begin{align}
	& 
	\a, \b > 0.
	\label{ab}
	\\ &
	\bz, \bQ, b_2, b_3   \, \hbox{ are nonnegative constants, but not all zero}.
	\label{constants}
	\\ &
	\phQ, \sQ \in L^2(Q), \sO \in \Huno, u_*,u^*\in L^\infty(Q) \hbox{ with } u_*\leq u^* \,\aeQ.
	\label{targets}
	\\ &
 	P \in C^2(\erre) \, \hbox{ is nonnegative, bounded and Lipschitz continuous}.
	\label{P}
	\\ &
	 \ph_0 \in \Hx3 \cap W, \m_0 \in \Huno, \s_0 \in \Huno.
	\label{initial_data}
\end{align}
Furthermore, we employ the following notation
\Bsist
	&& \non
	\UR \subset L^2(Q) \hbox{ be a non-empty and bounded open set such that
	it contains } \Uad
	\\ \non &&
	 \hbox{and } \norma u_2 \leq R \hbox{ for all }
	u \in \UR.
\Esist
For the potential setting, we require that $D(\hat{B})=\erre$ and that
\begin{align}
	& 
	\hat{B}:\erre \to [0, +\infty) \quad \hbox{is convex and lower semicontinuous, }
	\label{Bhat}
	\hbox{with $0 \in B(0)$}.
	\\ & \non
	\hat{\pi} \in C^1(\erre) \,\, \hbox{is nonnegative,} \,\, \pi:=\hat{\pi}' 
	\hbox{ is \Lip~continuous with}
	\\ &\label{pihat}
	\quad \hbox{\Lip~constant $L$, i.e. $\norma{\pi'}_{L^\infty(\erre)}\leq L$.}
\end{align}
Then, we define the potential $F$, and its derivative as the sum of these two contributions by
\Beq
	\label{defF}
	F:\erre \to [0, +\infty],  \quad  
	F:= \hat{B} + \hat{\pi} \aand F':= B + \pi,
\Eeq
where $B$ is a maximal and monotone graph $B\subset \erre \times \erre$ defined 
as the subdifferential of $\hat{B}$, that is, $B:=\partial\hat{B}$.
Unfortunately, we are not able to face the asymptotic analyses, as $\b$ 
goes to zero, without assuming proper growth restrictions for the potential $F$.
Some sufficient conditions for our purposes are as follows
\begin{align}
	\label{regpot}
	F &=\hat{B}+\hat{\pi} \quad \hbox{is a $C^3$ function which satisfies}
	\\	
	\label{growth}
	|B(r)| &\leq C_B (1 + \hat{B}(r)) \quad  \hbox{for every $r \in \erre$,}
\end{align}
\Accorpa\crescita regpot growth
\Accorpa\ipotesiF Bhat growth
\Accorpa\tutteleipotesi ab growth
for a positive constant $C_B$. 
Anyhow, we emphasize that, although we cannot work at the utmost
generality for the potentials setting, 
all polynomially growing potentials, as well as exponential functions, 
comply with the requirements above; in particular, \eqref{F_reg} is allowed.
Furthermore, by combining the embedding $W\subset \Lx\infty$ with the
first of the initial conditions \eqref{initial_data},
it is \sfw~to infer that $F(\ph_0)$ belongs to $\Lx\infty$.
It also follows from the above framework
that $F''$ is bounded below in terms of the \Lip\ constant $L$.
Indeed, we have that
\Beq
	\label{Fsecondoczero}
	F'' \geq - L.
\Eeq
It is worth noting that, in the case of \eqref{F_reg}, 
we can take $L=1$, as can be easily checked
by computing its second derivative.

Now, we first recall some results already presented in other contributions
and then list our statements. 
The already mentioned optimal control problem $(CP)_\b$ has been
tackled in \cite{S}. On the other hand,
since the above setting perfectly fits with
the one of \cite{S}, all the results there proved are at our disposal.
There, the author, after showing the existence
of optimal controls, provides the first-order necessary condition for optimality
making use of the \socal~adjoint system to \EQab. 
For the sake of simplicity, we just recall here the adjoint 
system there founded. From now on, $\uopt_\b$ denotes an optimal control
for $(CP)_\b$, and
$(\bm_\b,\bph_\b,\bs_\b)$ the corresponding optimal state. Thus, the adjoint system
reads as follows
\Bsist
  & \non
  \b \dt q_{\b} - \dt p_{\b} + \Delta q_{\b} - F''(\bph_{\b})q_{\b} 
  		+ P'(\bph_{\b})(\bs_{\b} - \bm_{\b})(r_{\b}-p_{\b})
  \\ \shoveright{ & = \bQ(\bph_{\b} - \phQ) \quad \hbox{in $\, Q$} }
  \label{EQAggabprima}
  \\
  & q_{\b} -\a \dt p_{\b} - \Delta p_{\b} + P(\bph_{\b})(p_{\b} - r_{\b}) = 0
  \label{EQAggabgseconda}
  \quad \hbox{in $\,Q$}
  \\
  & -\dt r_{\b} - \Delta r_{\b} + P(\bph_{\b})(r_{\b} - p_{\b}) = {\bQh (\bs_{\b} - \sQ)}
  \label{EQAggabterza}
  \quad \hbox{in $\,Q$}
  \\
  &\dn q_{\b} = \dn p_{\b} = \dn r_{\b} = 0
  \quad \hbox{on $\Sigma$}
  \label{BCEQAggab}
  \\ 
  & p_{\b}(T) - \b q_{\b}(T) =0,  \ \
  \a p_{\b}(T) = 0,  \  \
  r_{\b}(T) = {\bOh (\bs_{\b}(T) - \sO)}
  \quad \hbox{in $\Omega$.}
  \label{ICEQAggab}
\Esist
\Accorpa\EQAggab EQAggabprima ICEQAggab
The well-posedness of the above system has been already established
in \cite{S}. So, we just recall the obtained result.
\dafare{
\Bprop
Assume that \tutteleipotesi\ are in force. Then, system \EQAggab\
admits a unique solution satisfying 
\Beq
	\label{regagg}
	\qab, \pab, \rab \in  \H1 H \cap \L\infty V \cap \L2 W.
\Eeq
\Eprop}
Moreover, accounting for the solution of the adjoint system,
the following necessary condition was pointed out.
\Bthm
\label{THMsecondanecab}
Assume that \tutteleipotesi~are fulfilled. 
Let $\overline{u}_{\b} \in \Uad$ be an optimal control for $(CP)_{\b}$,
$(\bm_{\b},\bph_{\b},\bs_{\b})$ be the corresponding optimal state 
and $(p_{\b},q_{\b},r_{\b})$ the associated solution to the
adjoint system \EQAggab. Then, the necessary condition
for optimality is given by
\Beq
  \label{secondanecab}
  \intQ (r_{\b} + \bz \overline{u}_{\b})(v - \overline{u}_{\b})
  \geq 0 \quad \forall v \in \Uad.
\Eeq
\Ethm
As sketched above, we would like to exploit the control problem
$(CP)_\b$ in order to solve $(CP)$ by employing 
some asymptotic schemes.
In fact, in Section~\ref{SEC_CONTROL_PROBLEM},
we rigorously show that, as $\b\searrow 0$,
system \EQAggab~converge to
\Bsist
  & 
  	- \dt p 
    + \Delta q 
    - F''(\bph)q 
  	+ P'(\bph)(\bs - \bm)(r-p)
	= \bQ(\bph - \phQ) 
	\quad \hbox{in $\, Q$}
  \label{EQAggaprima}
  \\
  & q -\a \dt p - \Delta p + P(\bph)(p - r) = 0
  \label{EQAggagseconda}
  \quad \hbox{in $\,Q$}
  \\
  & -\dt r - \Delta r + P(\bph)(r - p) = {\bQh (\bs - \sQ)}
  \label{EQAggaterza}
  \quad \hbox{in $\,Q$}
  \\
  &\dn q = \dn p = \dn r = 0
  \quad \hbox{on $\Sigma$}
  \label{BCEQAgga}
  \\
  &  \a p(T) = 0,   \ \ r(T) = {\bOh (\bs(T) - \sO)}
  \quad \hbox{in $\Omega$.}
  \label{ICEQAgga}
\Esist
\Accorpa\EQAgga EQAggaprima ICEQAgga
Below, you can find the precise meaning of the above sentence.
\Bthm
\label{THM_stimeaggiuntoab}
Assume that \tutteleipotesi~are fulfilled and let 
$(\qab,\pab,\rab)$ be the unique solution to system \EQAggab~satisfying \eqref{regagg}. 
Then, \dafare{there exists $\a_{0}\in(0,1)$ such that, for every $\a\in(0,\a_{0})$},
as $\b \searrow 0$, and up to a subsequence, we have the following
convergences
\begin{align}
	\label{q_abtoqa}
	\qab &\to q \ \ \hbox{weakly in } \L2 W
	\\ 
	\label{p_abtopa}
	\pab &\to p \ \ \hbox{weakly star in } \H1 H \cap \L\infty V \cap \L2 W
	\\   \label{r_abtora}
	\rab &\to r \ \ \hbox{weakly star in } \H1 H \cap \L\infty V \cap \L2 W
	\\   \label{bqtozero_due}
	\b q_{\b} &\to 0 \ \ \hbox{strongly in } \H1 H \cap \L\infty V \cap \L2 W.
\end{align}
\Accorpa\convadjab q_abtoqa bqtozero_due
Moreover, there exists a positive constant $C_1$, independent of $\b$, such that
\begin{align}
	\label{stimaagg}
	& \non
	\b\norma {\qab}_{\H1 H}
	+ \b^{1/2} \norma {\qab}_{\L\infty V}
	+ \norma{\qab}_{\L2 W}
	\\ & \quad 
	+ \norma{\pab}_{\H1 H \cap \L\infty V\cap \L2 W}
	+ \norma{\rab}_{\H1 H \cap \L\infty V \cap \L2 W}
	\leq C_1.
\end{align}
In addition, the limit $(q,p,r)$ is the unique solution to system
\EQAgga~which possesses the following regularity 
\begin{align}
	\label{reg_adj_lim_uno}
	q &\in \L2 W
	\\ 
	p, r &\in \H1 H \cap \L\infty V \cap \L2 W.
	\label{reg_adj_lim_due}
\end{align}
\Accorpa\regadjlim {reg_adj_lim_uno} {reg_adj_lim_due}
\Ethm
Next, we can address the results related to the control problem 
we are dealing with.
We begin with the first fundamental result
concerning the existence of optimal control.
\Bthm
\label{THexistenceofcontrol}
Assume that \tutteleipotesi~are in force. 
Then, the optimal control problem $(CP)$ admits at least a solution
$\overline{u} \in \Uad$.
\Ethm
Lastly, by employing a proper asymptotic scheme, we develop the first-order
necessary condition for optimality.
\Bthm
\label{THMsecondanec}
Assume that \tutteleipotesi~are satisfied. 
Let $\overline{u} \in \Uad$ be an optimal control for $(CP)$
with its corresponding state $(\bm,\bph,\bs)$ and 
let $(p,q,r)$ be the solution to the associated adjoint system
\EQAgga. Then, the following variational inequality
\Beq
  \label{secondanec}
  \intQ (r + \bz \overline{u})(v - \overline{u})
  \geq 0 \quad \forall v \in \Uad
\Eeq
is satisfied.
Moreover, whenever $\bz \not= 0$, 
the optimal control $\overline{u}$ is the $\L2 H-$projection
of $-{r}/{\bz}$ onto the subspace $\Uad$.
\Ethm

To conclude the section, we recall a \wk~inequality and a general fact
that is widely used in the sequel.
First of all, let us remind the Young inequality
\Beq
  ab \leq \delta a^2 + \frac 1 {4\delta} \, b^2
  \quad \hbox{for every $a,b\geq 0$ and $\delta>0$}.
  \non
  \label{young}
\Eeq
Furthermore, we recall the standard Sobolev continuous embedding
\Beq
	\label{VinL6}
	\Huno \hookrightarrow L^q(\Omega)
	\quad \hbox{which holds for every $q\in [1,6]$.}
\Eeq
Throughout the paper, we convey to use the symbol small-case
$c$ for every constant
which depend only on the final time~$T$, on~$\Omega$, on $R$, on the shape of the 
nonlinearities, on the norms of the involved functions, and possibly on $\a$. 
On the other hand, we will explicitly point out when an appearing 
constant may depend on $\b$.
For this reason, the meaning of $c$ might
change from line to line and even in the same chain of inequalities.
Differently, we devote the capital letters
to indicate particular constants which we eventually will refer later on.
\section{The State System}
\label{SEC_STATE}
\setcounter{equation}{0}

From this section onward, we start with the proofs of the introduced
statements. We begin with the task of investigating the well-posedness of system \EQab,
and proving its asymptotic behavior as $\b\searrow 0$. 
Again, we remark that the above system has already been
discussed in \cite{CGH, CGRS_ASY, CGRS_VAN}, where the asymptotic analysis
represents the core of the contributions. For this reason, some of the calculations
below can also be found there.
Anyhow, in order to handle the control problem $(CP)$, 
it turns out that the results there obtained are insufficient. 
Therefore, we perform all the necessary estimates, 
having the care to emphasize when 
the appearing constants may depend on $\b$.

\Bthm
\label{THM_exist_ab}
Let the assumptions \tutteleipotesi~be fulfilled \dafare{and let $\m_0:=-\Delta \ph_0 +F'(\ph_0)$.}
Then, \dafare{there exists $\a_{0}\in(0,1)$ such that, for every $\a\in(0,\a_{0})$},
system \EQab~admits a unique solution $(\m_{\b}, \ph_{\b}, \s_{\b})$
that satisfies the following regularity
\begin{align}
	\label{reg_thm_uno}
	\m_\b &\in \H1 H \cap \L\infty  V \cap \L2 W
	\\  
	\label{reg_thm_due} 
	\ph_\b &\in W^{1,\infty}(0,T; H) \cap \H1 V \cap \L\infty W
	\\   
	\label{reg_thm_tre}
	\s_\b& \in	 \H1 H \cap \L\infty  V \cap \L2 W.
\end{align}
\Accorpa\regthm  reg_thm_uno reg_thm_tre
Furthermore, there exists a positive constant $C_2$
such that the following estimate is verified
\begin{align}
	\non &
	\b^{1/2}\norma{\dt \ph_{\b}}_{\L\infty H}
	+ \norma{\ph_{\b}}_{ \H1 V \cap \L\infty W \cap L^\infty( Q)} 
	\\ &\qquad \non 
	+ \norma{\m_{\b}}_{\H1 H \cap \L\infty V \cap \L2 W}
	+ \norma{\s_{\b}}_{\H1 H \cap \L\infty V \cap \L2 W}
	\\ &\quad 
	\leq
	{C_2} \bigl(
	\norma{\m_{0}}_{V}
	+ \norma{\ph_{0}}_{W}
	+ \norma{\s_{0}}_{V}
	+1
	\bigr),
	\label{stima_ab_uno}
\end{align}
where ${C_2}$ is a positive constant that depends on $\Omega, T, R, \a$,
the shape of the functions $P$ and $\pi$, but it is independent of $\b$.
\Ethm

Before moving to the proof,
it is worth focusing the attention on the already obtained
limit system \EQa~which was investigated in \cite[Thm.~2.2, p.~41]{CGRS_ASY}
(see also \cite{CGRS_VAN }).
Forgetting the fact that \EQa~has founded as a result of a limit procedure, 
it can be considered as a starting point as a system of partial differential 
equations. 
\dafare{In this regards, we only have at disposal the derivative 
of the linear combination $\a\m+\ph$, so that the initial condition, in general, reads as 
$(\a\m+\ph) (0)=\et_0$,
}
for a suitable element $\et_0$.
Furthermore, owing to what already proved for the 
system \EQa, we claim that, whenever $\et_0\in V$,
the existence and the uniqueness of the solution follows, provide
a smallness assumption of the constant $\a$ is satisfied. 
Thus, we consider this problem under
the additional assumptions
\Beq
\dafare{\begin{cases}
	 \et_0 = \a \m_0 + \ph_0 
     \\  
     \m_0  = - \Delta \ph_0 + F'(\ph_0).
     \label{def_m0}
\end{cases}}
\Eeq
The first condition
allows us to match the two approaches; that is, to see the above system as
a direct problem and as the limit of system \EQab~as $\b\searrow 0$.
Conversely, the second condition states that
$\m_{0}$ cannot be arbitrary chosen, but it has to be defined in terms of $\ph_{0}$
in a prescribed way. At this stage, that requirement 
could appear quite unnecessary and unnatural, but it will be motivated in 
view of a forthcoming estimate (see the Fourth estimate below)
which is of crucial importance for the asymptotic analysis. 
Under this additional assumption, it follows that
the second of \eqref{initial_data} is rather a consequence of the first
one combined with the strong regularity we postulated for $\ph_0$, which 
can now be motivated by virtue of the regularity we want for $\m_{0}$. 
Indeed, combining the growth assumption for the potential \eqref{growth} with 
the continuous embedding $W\subset L^\infty(\Omega)$,
we infer from the second of \eqref{def_m0} that $\m_0\in V$.
On the other hand, whenever $\et_0\in V$ is given, 
let us claim that $\m_0$ and $\ph_0$ can be reconstructed providing
to impose $\dn \ph_0=0$. In fact, using \eqref{def_m0}, we realize that
we are looking for a variable $\ph_0$ such that 
\Beq
	\non
	\begin{cases}
	- \Delta \ph_0 + \frac{\ph_0-\et_0}{\a} + F'(\ph_0) = 0 \quad \hbox{in $\Omega$}
	\\
	\dn \ph_0=0	\quad \hbox{on $\Gamma$}.
	\end{cases}
\Eeq
Again, the fact that $\a$ has to be sufficiently small helps us and 
the existence and uniqueness of a solution to the above equation can be proved.
Indeed, the nonlinear term $F'(\ph_0)$ can be split as 
$F'(\ph_0)= B(\ph_0)+\pi(\ph_0)$, where 
the perturbation $\pi(\ph_0)$ does not bother since it is balanced by the term
$\frac{\ph_0-\et_0}{\a}$ which dominates owing to the smallness of the denominator
(we can indeed assume that $\a L < 1$, which gives $\frac 1\a > L$). 
Then, by combining $\partial_n \ph_0=0$ with the elliptic regularity theory, 
we are able to reconstruct $\ph_0$ which fulfills the first condition 
of \eqref{initial_data}. Finally, from the first equation of \eqref{def_m0}, we also recover $\m_0$
with the prescribed regularity.
With these comments 
in mind, we can proceed with the proof of Theorem \ref{THM_exist_ab}.

\proof[Proof of Theorem \ref{THM_exist_ab}]
For the uniqueness part, we refer the reader to \cite[Sec.~3]{CGH}.
As the existence is concerned, we take into account an approximation scheme.

{\bf The approximating system}
Let us take $\eps \in (0,1)$ and consider the \socal~Yosida approximation of the 
maximal and monotone operator ${B}$. Namely, for every $r \in \erre$, we introduce
\Bsist
	\label{yosida}
	\hat{B}_\eps (r) := \min_{s\in\erre}  \biggl( \frac 1{2\eps} (s-r)^2 + \hat{B}(s) \biggr),	
	\quad B_{\eps} (r) := \frac d{dr} \hat{B}_\eps (r), \aand
	F_\eps:= \hat{B}_\eps  + \hat{\pi}.
\Esist
It turns out that $\hat{B}_\eps$ is a well-defined $C^1$ function,
$B_\eps$ is \Lip~continuous (see, e.g., \cite[Prop.~2.11, p.~39]{BRZ}),
and for every $r\in\erre$, it holds that
\Bsist
	\label{yosida_mono}
	0 \leq \hat{B}_\eps (r) \leq \hat{B} (r) \quad \hbox{and} \quad
	\hat{B}_\eps (r) \nearrow \hat{B} (r) \quad  \hbox{as $\eps\searrow 0$.}
\Esist
Hence, in order to solve \EQab, we first solve the approximated
system obtained by substituting $F$ by $F_\eps$, and then we 
let $\eps\searrow 0$ to prove the existence to the original
problem. Therefore, we are going to face the following system
\Bsist
  & \a \dt \m_{\b,\eps} + \dt \ph_{\b,\eps} 
  	- \Delta \m_{\b,\eps} = P(\phi_{\b,\eps}) (\s_{\b,\eps} - \m_{\b,\eps})
  \quad \hbox{in $\, Q$}
  \label{EQAPPprima}
  \\
  & \m_{\b,\eps} = \beta \dt \ph_{\b,\eps} - \Delta \ph_{\b,\eps} + F_{\eps}'(\ph_{\b,\eps})
  \label{EQAPPseconda}
  \quad \hbox{in $\,Q$}
  \\
  & \dt \s_{\b,\eps} - \Delta \s_{\b,\eps} = - P(\ph_{\b,\eps}) (\s_{\b,\eps} - \m_{\b,\eps}) + u_\b
  \label{EQAPPterza}
  \quad \hbox{in $\,Q$}
  \\
  & \dn \m_{\b,\eps} = \dn \ph_{\b,\eps} =\dn \s_{\b,\eps} = 0
  \quad \hbox{on $\,\Sigma$}
  \label{BCEQAPP}
  \\
  & \m_{\b,\eps}(0)=\m_{0},\, \ph_{\b,\eps}(0)=\ph_{0},\, \s_{\b,\eps}(0)=\s_{0}
  \quad \hbox{in $\,\Omega.$}
  \label{ICEQAPP}
\Esist
\Accorpa\EQAPP EQAPPprima ICEQAPP

Our starting point is the result below.
\Blem
Assume that \tutteleipotesi~are satisfied. Then, the approximating 
problem \EQAPP~admits a unique solution.
\Elem
As the uniqueness is concerned, it can be proved as a special case of
\dafare{\cite[Sec.~3]{CGH}. As regards the existence,}
let us only mention that, e.g., a Faedo-Galerkin scheme,
along with some a priori estimates, will lead to proving the asserted result.
For instance, as a basis of $V$, one could take into account
the basis consisting of the eigenfunctions of the Laplacian operator
with homogeneous Neumann boundary conditions.
We decide to skip the details because the estimates we are
going to perform below are very similar to the ones that 
could allow one to solve the approximating problem.
In addition, let us point out that ${B}_\eps$ and the map which assigns
$(\m_{\b,\eps},\ph_{\b,\eps},\s_{\b,\eps})
\mapsto P(\ph_{\b,\eps})(\m_{\b,\eps}-\s_{\b,\eps})=:R_{\b,\eps}$
are both smooth and \Lip~continuous, for every $\b$ and every $\eps$, 
and therefore the classical
Picard-Lindel\"of theorem directly yields the existence of a unique global solution
to the system of ordinary differential equations given by that scheme.

{\bf First estimate} 
We multiply \eqref{EQAPPprima} by $\m_{\b,\eps}$, \eqref{EQAPPseconda} by 
$-\dt\ph_{\b,\eps}$, and \eqref{EQAPPterza}
by $\s_{\b,\eps}$. Next, we integrate over $Q_t$ and by parts, and add the 
resulting equalities to obtain that
\Bsist
	\non &&
	\frac \a2 \IO2 {\m_{\b,\eps}}
	+ \I2 {\nabla\m_{\b,\eps}}
	+ \b \I2 {\dt\ph_{\b,\eps}}
	+ \frac 12 \IO2 {\nabla\ph_{\b,\eps}}
	\\ \non && \qquad
	+ \iO {F}_\eps(\ph_{\b,\eps}(t))
	+ \frac 12 \IO2 {\s_{\b,\eps}}
	+ \I2 {\nabla\s_{\b,\eps}}
	+ \intQt P(\ph_{\b,\eps})(\s_{\b,\eps}-\m_{\b,\eps})^2
	\\ \non && \quad
	= 
	\frac \a2 \iO |\m_{0}|^2
	+ \frac 12 \iO |\nabla\ph_{0}|^2
	+ \iO {F}_\eps(\ph_{0})
	+ \frac 12 \iO |\s_{0}|^2
	+ \intQt u_\b \s_{\b,\eps},
\Esist
where we denote the integrals on the \rhs~by $I_1,...,I_5$, in this order.
The terms on the \lhs~are nonnegative since they all are squares and  
$P$ and $F_\eps$ are so by \eqref{P} and \eqref{pihat} 
along with \accorpa{yosida}{yosida_mono}, respectively.
Moreover, $I_1,I_2$ and $I_4$ can be easily controlled owing to the 
assumptions on the initial conditions \eqref{initial_data}. 
As for $I_3$, we deduce that
\Bsist
	\non
	|I_3|
	=
	\Big|\iO {F}_\eps(\ph_{0})\, \Big|
	= 
	\iO {\hat{B}_\eps}(\ph_{0})
	+
	\iO {\hat{\pi}(\ph_{0})}
	\leq 	
	\iO {\hat{B}}(\ph_{0})
	+
	\iO {\hat{\pi}(\ph_{0})}
	\leq 
	c,
\Esist
by invoking the properties of $\hat{B}_\eps$ pointed out by \eqref{yosida_mono} 
and accounting for the 
features of the initial datum $\ph_{0}$ and on the function $\hat{\pi}$.
Then, by employing the Young inequality, we bound $I_5$ as follows
\Bsist
	\non
	|I_5|
	\leq
	\frac 12 \I2 {u_\b}
	+ \frac 12 \I2 {\s_{\b,\eps}}.
\Esist
Thus, a Gronwall argument yields that
\begin{align}
 	\non &
 	\norma{\m_{\b,\eps}}_{\L\infty H \cap \L2 V}
	+ \b^{1/2}\norma{\dt\ph_{\b,\eps}}_{\L2 H }
	+ \norma{\nabla \ph_{\b,\eps}}_{\L2 H }
	\\  & \quad
	+ \norma{F_\eps(\ph_{\b,\eps})}_{\L\infty \Luno}
	+ \norma{\s_{\b,\eps}}_{\L\infty H \cap \L2 V}
	\leq 
	c.
	\label{Iest}
\end{align}

{\bf Second estimate}
Analyzing \eqref{Iest}, we realize that it does not
provide any information of $\ph_{\b,\eps} $ in $\L2 H$. 
Hence, we try to reconstruct the full norm of $\ph_{\b,\eps}$ in $\L2 V$. 
In this direction, we add equations \eqref{EQAPPprima} and \eqref{EQAPPterza}
to get 
\Beq
	\dt (\a \m_{\b,\eps} + \ph_{\b,\eps} + \s_{\b,\eps} )
	- \Delta (\m_{\b,\eps} + \s_{\b,\eps}) = u_\b.
	\non
\Eeq
Then, we test the above equation by $\a \m_{\b,\eps} + \ph_{\b,\eps} + \s_{\b,\eps}$
and integrate over $Q_t$ and by parts. 
Upon rearrange the terms, we obtain that
\begin{align}
	& \non
	\frac 12 \iO (\a |\m_{\b,\eps} (t)|^2 + |\ph_{\b,\eps} (t)|^2 + |\s_{\b,\eps} (t)|^2 )
	+ \a \I2 {\nabla\m_{\b,\eps}}
	+  \I2 {\nabla\s_{\b,\eps}}
	\\ & \quad \non
	=
	\frac 12 \iO (\a |\m_{0}|^2 + |\ph_{0}|^2 + |\s_{0} |^2 )
	- (\a+1) \intQt \nabla\m_{\b,\eps} \cdot \nabla \s_{\b,\eps}
	\\ & \qquad \non
	- \intQt \nabla\m_{\b,\eps} \cdot \nabla \ph_{\b,\eps}
	- \intQt \nabla\s_{\b,\eps} \cdot \nabla \ph_{\b,\eps}
	+ \intQt u_\b (\a \m_{\b,\eps} + \ph_{\b,\eps} + \s_{\b,\eps} ),
\end{align}
where the integrals on the \rhs~are denoted by $I_1,\dots, I_5$, in this order.
Using \eqref{initial_data}, we immediately deduce that $|I_1|\leq c$.
Furthermore, by combining the above estimate with the Young inequality, we find that the 
remaining terms can be estimated as
\Bsist
	&& |I_2|+|I_3|+|I_4|+|I_5|
	\leq \non
	c \intQt (|\m_{\b,\eps}|^2 + |\nabla \m_{\b,\eps}|^2)
	+ c \intQt (|\ph_{\b,\eps}|^2 + |\nabla \ph_{\b,\eps}|^2)
	\\ && \quad \non
	+ c \intQt (|\s_{\b,\eps}|^2 + |\nabla \s_{\b,\eps}|^2)
	+ \frac 12 \I2 {u_\b}.
\Esist
Therefore, the Gronwall lemma entails that
\Bsist
 	&&
 	\norma{\ph_{\b,\eps}}_{\L\infty H }
	\leq 
	c.
	\label{IIest_a}
\Esist
Moreover, we also realize that
\Beq
	\label{IIest_b}
	\norma{R_{\b,\eps}}_{\L2 H }
	\leq 
	c.
\Eeq
\Accorpa\IIest IIest_a IIest_b

{\bf Third estimate}
It is worth noting that \eqref{EQAPPterza} possesses a parabolic structure with respect to the
variable $\s_{\b,\eps}$ and,
owing to the above results, it follows that its forcing term belongs to $\L2 H$.
Hence, the parabolic regularity theory for homogeneous Neumann problems 
with regular initial conditions, gives us
\Bsist
	\label{IIIest}
	\norma{\s_{\b,\eps}}_{\H1 H \cap \L\infty V \cap \L2 W}
	\leq
	c.
\Esist

{\bf Fourth estimate}
Now, we present the key estimate for the forthcoming asymptotic analyses
which motivates the unusual requirement \eqref{def_m0}.
To begin with, let us formally differentiate \eqref{EQAPPseconda} with
respect to time to infer that
\Bsist
	\label{dtseconda}
	\dt\m_{\b,\eps} = \beta \dtt \ph_{\b,\eps} - 
	\Delta \dt\ph_{\b,\eps} + F_{\eps}''(\ph_{\b,\eps})\dt\ph_{\b,\eps}.
\Esist
Next, we multiply it by $\a$ and replace the first term of \eqref{EQAPPprima}
with this new equation leading to obtain that
\begin{align}
	\label{eqformale}
	\a \b \dtt \ph_{\b,\eps} - 
	\a\Delta \dt\ph_{\b,\eps} + \a F_{\eps}''(\ph_{\b,\eps})\dt\ph_{\b,\eps}
	+ \dt \ph_{\b,\eps} - \Delta \m_{\b,\eps} = P(\phi_{\b,\eps}) (\s_{\b,\eps} - \m_{\b,\eps}).
\end{align}
This formal procedure can be rigorously motivated.
Indeed, by introducing the auxiliary variable $z_\b:=\a \dt\ph_{\b,\eps}$, 
we can rewrite \eqref{eqformale} as a 
parabolic equation as follows
\Beq
	\non
	\b\dt z_\b - \Delta z_\b = f_\b \quad \hbox{in $Q$}, 
\Eeq
where $f_\b$ is defined by
\Beq
	\non
	f_\b:= 
	\Delta\m_{\b,\eps}
	- \dt\ph_{\b,\eps}
	- \a F_{\eps}''(\ph_{\b,\eps}) \dt\ph_{\b,\eps}
	+  P(\phi_{\b,\eps}) (\s_{\b,\eps} - \m_{\b,\eps}).
\Eeq
Owing to the above estimates and to the growth conditions \crescita~for the potential,
we easily realize that, for every $\b$, the forcing term $f_\b \in \L2 \Vp$. 
Therefore, the abstract theory for parabolic equations
(see, e.g., \cite{Lions})
guarantees the existence and the uniqueness of a solution $z_\b\in \H1 \Vp \cap \L2 V$, whenever
the initial datum $z_\b(0)$ is sufficiently regular, that is, whenever $z_\b(0)$
belongs at least to $H$.
As we will see, the particular choice of the initial 
datum $\m_0$ assumed by \eqref{def_m0}, entails, in the limit,
that $z_\b(0)=0$, so that the required
regularity is trivially fulfilled.

Then, we multiply \eqref{eqformale} by $\dt\ph_{\b,\eps}$ and integrate over $Q_t$ and 
by parts to find that
\begin{align}
	\non &
	\frac{\a \b}2  \IO2 {\dt\ph_{\b,\eps}}
	+ \a \I2 {\nabla \dt\ph_{\b,\eps} }
	+ \a \intQt F_{\eps}''(\ph_{\b,\eps}) |{\dt\ph_{\b,\eps} }|^2
	+ \intQt |{\dt\ph_{\b,\eps} }|^2	
	\\ & \quad 
	= 
	\frac{\a \b}2  \iO |\dt\ph_{\b,\eps} (0)|^2
	+ \intQt P(\phi_{\b,\eps}) (\s_{\b,\eps} - \m_{\b,\eps})\,\dt\ph_{\b,\eps}
	- \intQt \nabla \m_{\b,\eps} \cdot \nabla \dt\ph_{\b,\eps},
	\label{system_in}
\end{align}
where the integrals on the \rhs~are denoted by $I_1,I_2,I_3$, in this order.
As the third integrals of the \lhs~is concerned,
we remind that $F_{\eps}$ is defined by \eqref{yosida} and that $F''$
is bounded below. 
Moreover, the same bound holds true for $F_\eps$ 
with the same constant $L$. Therefore,
the third and fourth contributions on the \lhs~verify that
\Beq
	\label{aFdue}
	\a \intQt F_{\eps}''(\ph_{\b,\eps}) |{\dt\ph_{\b,\eps} }|^2
	+ \intQt |{\dt\ph_{\b,\eps} }|^2	
	\geq 
	(1-\a L)	\intQt |{\dt\ph_{\b,\eps} }|^2,
\Eeq
whereas the other terms on that side are nonnegative.
As regards the \rhs, let us emphasize that the definition \eqref{def_m0} 
plays a fundamental role.
Indeed, by taking $t=0$ in the equation \eqref{EQAPPseconda}, we get
\Beq
	\non
	 \m_{\b,\eps}(0) = \beta \dt \ph_{\b,\eps}(0) - 
	 \Delta \ph_{\b,\eps}(0) + F_{\eps}'(\ph_{\b,\eps}(0))
	 \quad
	 \hbox{in $\Omega$}
\Eeq
Using \eqref{def_m0}, the initial conditions \eqref{ICEQAPP}, and
rearranging the terms lead to infer that
\Beq
	\non
	\beta \dt \ph_{\b,\eps}(0) 
	= \m_0 + \Delta \ph_{0} - F_{\eps}'(\ph_0)
	= B(\ph_0)-B_{\eps}(\ph_0)
	\quad
	\hbox{in $\Omega$.}
\Eeq
\dafare{
On the other hand, from \wk~results (see, e.g., \cite[Chapt.~2]{BRZ}),
we also have the convergence, both $a.e.$ in $\Omega$ and in $\Lx2$, of
the operator $B_{\eps}$ to element having
minimum norm of the limit operator $B$. Actually, since $\hat{B}$ is regular, 
it turns out that its subdifferential $B$ is single-valued, so that the element having
minimum norm is $B$ itself. Hence, we have
\Beq
	B_{\eps}(r) \to B(r) \quad \hbox{strongly in $\Lx2$, $a.e.$ in $\Omega$, for every $r \in \erre$.}
\Eeq
Therefore, we easily realize that $\normaH{B(\ph_0)-B_{\eps}(\ph_0)}^2\to 0$, 
as $\eps \to 0$.}
Meanwhile, the other integrals can be easily controlled thanks to the Young
inequality. Recalling that $P$ is bounded by \eqref{P}, we control $I_2$ by
\Beq
	\label{I2}
	|I_2|
	\leq
	\d \I2 {\dt\ph_{\b,\eps}}
	+ \cd  \intQt ({|\s_{\b,\eps} |^2 + |\m_{\b,\eps} |^2}),
\Eeq
for a positive constant $\d$, yet to be determined. 
In a similar manner, we obtain
\Beq
	\label{I3}
	|I_3|
	\leq
	\frac \a 2  \intQt |\nabla \dt\ph_{\b,\eps}|^2
	+ c \I2 {\nabla\m_{\b,\eps}}.
\Eeq
Upon collecting \accorpa{aFdue}{I3}, we rearrange the initial equation 
\eqref{system_in} to infer that
\begin{align}
	\non &
	\frac{\a \b}2  \IO2 {\dt\ph_{\b,\eps}}
	+ \frac \a 2 \I2 {\nabla \dt\ph_{\b,\eps} }
	+ (1-\a L - \d)\intQt |{\dt\ph_{\b,\eps} }|^2
	\\ \non 
	& \quad \leq 
	c(1+ \normaH{B(\ph_0)-B_{\eps}(\ph_0)}^2)
\end{align}
has been shown, 
where the \rhs~has been managed owing to the above estimates.
On the other hand, we can assume $\a$ to be sufficiently small 
in order to have that $\a L<1$, so that
it suffices to take $\d$ small enough to conclude $(1-\a L - \d)>0$, which, in turn,
implies that
\Bsist
	\label{IVest}
	\b^{1/2} \norma{\dt\ph_{\b,\eps}}_{\L\infty H}
	+ \norma{\dt\ph_{\b,\eps}}_{\L2 V}
	\leq 
	c(1+ \normaH{B(\ph_0)-B_{\eps}(\ph_0)}).
\Esist
%
{\bf Fifth estimate}
Moreover, \eqref{EQAPPprima} shows a parabolic structure with respect
to the variable $\m_{\b,\eps}$, so that
\Bsist
	\non
	\a \dt \m_{\b,\eps} - \Delta \m_{\b,\eps} 
	= f_\b  \quad \hbox{in $Q$,}
\quad \hbox{where } f_\b:= -\dt \ph_{\b,\eps} + P(\phi_{\b,\eps}) (\s_{\b,\eps} - \m_{\b,\eps}).
\Esist
On account of the above estimates, we easily realize that
$f_\b\in\L2 H$.
Thus, the regularity theory for homogeneous Neumann parabolic equations yields that
\Bsist
	\label{Vest}
    \norma{\m_{\b,\eps}}_{\H1 H \cap \L\infty V \cap \L2 W}
	\leq 
	c.
\Esist

{\bf Sixth estimate}
Furthermore, let us read \eqref{EQAPPseconda} as an elliptic equation 
with respect to the variable $\ph_{\b,\eps}$ as follows
\Beq
	\non
	- \Delta \ph_{\b,\eps} + F'_\eps(\ph_{\b,\eps}) = \m_{\b,\eps} - \b \dt \ph_{\b,\eps} 
	\quad
	\hbox{in $Q$}.
\Eeq
Then, we consider the above equation written at time $t$,
split $F'_\eps$ on account of \eqref{defF}, multiply it by
$-\Delta \ph_{\b,\eps}(t)$ and integrate over $\Omega$ and by parts to obtain that
\begin{align}
	& \non
	\IO2 {\Delta \ph_{\b,\eps}}
	+ \iO B'_\eps(\ph_{\b,\eps}(t)) |\nabla \ph_{\b,\eps}(t)|^2 
	=
	- \iO \m_{\b,\eps}(t) \, \Delta \ph_{\b,\eps}(t)
	\\ & \qquad \non
	+ \b \iO \dt \ph_{\b,\eps}(t) \, \Delta \ph_{\b,\eps}(t)
	+ \iO \pi( \ph_{\b,\eps}(t)) \, \Delta \ph_{\b,\eps}(t),
\end{align}
where the terms on the \rhs~are denoted by $I_1,I_2$
and $I_3$, in that order.
Note that, at the first stage, 
the second term on the \rhs~can be neglected since it is nonnegative
by the properties of $B'_\eps$. On the other hand, 
Young's inequality, along with the above estimates, gives
\Beq
	\non
	|I_1|+|I_2|+|I_3|
	\leq
	\frac 34 \iO |\Delta \ph_{\b,\eps}(t)|^2
	+c.
\Eeq
Hence, invoking first the elliptic theory, and then comparison
in \eqref{EQAPPseconda}, lead to conclude that
\Beq
	\label{VIIest_a}
    \norma{\ph_{\b,\eps}}_{\L\infty W}
    + 
	\norma{B_\eps(\ph_{\b,\eps})}_{\L\infty H} 
	\leq 
	c.
\Eeq
Furthermore, the continuous embedding $W \subset L^\infty(\Omega)$, entails that
\Beq
	\label{VIIest_b}
    \norma{\ph_{\b,\eps}}_{L^\infty( Q)}
	\leq 
	c.
\Eeq
\Accorpa\VIIest VIIest_a VIIest_b

{\bf Passing to the limit}
Lastly, we draw some consequences from the above a priori estimates,
showing that we can let $\eps\searrow 0$
to complete the proof.

Owing to standard weak compactness arguments,
we infer that, as $\eps\searrow 0$, and up to a not relabeled subsequence,
the following convergences
\begin{align}
	\m_{\b,\eps} &\to \m_{\b} \ \ \hbox{weakly star in }  \H1 H \cap \L\infty V \cap \L2 W
	\\  
	\ph_{\b,\eps} &\to \ph_{\b}  \ \ \hbox{weakly star in } W^{1,\infty}(0,T; H) \cap \H1 V \cap \L\infty W
	\\  
	\s_{\b,\eps} &\to \s_{\b}  \ \ \hbox{weakly star in }  \H1 H \cap \L\infty V \cap \L2 W
	\\   
	R_{\b,\eps} &\to \zeta_{\b}  \ \ \hbox{weakly in } \L2 H
	\\  
	B_{\eps}(\ph_{\b,\eps}) &\to \psi_\b \ \ \hbox{weakly star in } \L\infty  H
\end{align}
are satisfied.
Furthermore, compactness embedding results (see, e.g., \cite[Sec.~8, Cor.~4]{Simon}) 
easily imply that
\begin{align}
	\non
	\m_{\b,\eps} \to \m_{\b}, \quad \ph_{\b,\eps} \to \ph_{\b}, \quad \s_{\b,\eps} \to \s_{\b}
	\quad \hbox{strongly in } \C0 H \cap \L2 V .
\end{align}
These convergences give sense to initial conditions \eqref{ICEQab} 
and allow us to identify the limit
of the nonlinear terms. Indeed, the boundedness and the 
regularity of $P$, along with the above strong convergences,
yield that $\zeta_{\b}=R_\b$, where $R_\b:=P(\ph_\b)(\s_\b-\m_\b)$. 
Arguing in a similar fashion, we infer that
$\pi(\ph_{\b,\eps})$ strongly converges to $\pi(\ph_{\b})$ in $\L2 H$.
Lastly, from the monotonicity properties of the Yosida approximation introduced by 
\eqref{yosida}, we get (see, e.g., \cite[Lemma~1.3, p.~42]{BAR}) that
\Beq
	\non
	\limsup_{\eps \searrow 0} \intQ B_{\eps}(\ph_{\b,\eps})\ph_{\b,\eps}
	=
	\lim_{\eps \searrow 0} \intQ B_{\eps}(\ph_{\b,\eps})\ph_{\b,\eps}
	=
	\intQ B(\ph_{\b})\ph_{\b},
\Eeq
which implies $\psi_\b=B(\ph_{\b}) $.
In conclusion, the limit triplet $(\m_\b,\ph_\b,\s_\b)$ yields a solution to \EQab~and 
possesses the postulated regularity \regthm, so that Theorem \ref{THM_exist_ab}
is completely proved.
\qed

With the result below, we aim at improve the result on the 
asymptotic analysis of system \EQab, as $\b \searrow 0$ 
(compare with the regularity pointed out in \cite{CGH,CGRS_VAN,CGRS_ASY}).
\Bthm
\label{THM_asy_dir}
Suppose that \tutteleipotesi~are satisfied \dafare{and let $\m_0:=-\Delta \ph_0 +F'(\ph_0)$.}
Then, there exists a sufficiently small $\a_{0}\in(0,1)$ such that, for every
$\a \in (0,\a_{0})$ and $\b\in (0,1)$, the unique
solution $(\m_{\b},\ph_{\b},\s_{\b})$ to problem 
\EQab, as $\b \searrow 0$, satisfies
\begin{align} 
	\label{mab_ma}
	\m_{\b} &\to \m \ \ \hbox{weakly star in } \H1 H \cap \L\infty V \cap \L2 W
	\\  
	\label{phab_pha}
	\ph_{\b} &\to \ph \ \ \hbox{weakly star in } {\H1 V \cap \L\infty W}
	\\   \label{sab_sa}
	\s_{\b} &\to \s \ \ \hbox{weakly star in } \H1 H \cap \L\infty V \cap \L2 W
	\\  
	\label{bphab_zero}
	\b \ph_{\b} &\to 0 \ \ \hbox{strongly in } {W^{1,\infty}(0,T; H) \cap \H1 V \cap \L\infty W}
\end{align}
at least for a subsequence. 
Moreover, the limit $(\m,\ph,\s)$ turns out to be the unique solution to the 
limit system \EQa. Furthermore, there exists a subsequence for which
we also have the strong convergences
\begin{align}
	\non
	\ph_{\b} & \to \ph \ \ \hbox{strongly in } \C0 {\Hx{2-\gamma}},\ \hbox{for every $\gamma > 0$,} 
	\\  
	\label{strong_phab}
	& \quad \hbox{which entails} \ \ \ph_{\b} \to \ph \ \ \hbox{strongly in } C^0 (\overline{Q})
	\\  
	\label{strong_mab}
	\m_{\b} &\to \m \ \ \hbox{strongly in } \C0 H\cap \L2 V 
	\\  
	\label{strong_sab}
	\s_{\b} &\to \s \ \ \hbox{strongly in } \C0 H\cap \L2 V .
\end{align}
\Ethm

\Bdim
Since it immediately follows on account of standard techniques from
estimate \eqref{stima_ab_uno},
we just sketch the proof and left the details to the reader.

The convergences \accorpa{mab_ma}{bphab_zero} immediately follow by standard argument
by estimate \eqref{stima_ab_uno}, while
the strong convergences \accorpa{strong_phab}{strong_sab}
can be obtained by invoking \wk~compactness results
(see, e.g., \cite{Simon}).
Hence, it remains to check that $(\m,\ph,\s)$ yields
a solution to \EQa. In principle, one should consider the variational formulation 
corresponding to system \EQab~and, using the above estimates, pass to the limit
to conclude. Therefore, the only terms that deserve further comments are the nonlinear ones.
Anyhow, the strong convergence \eqref{strong_phab} suffices since,
along with \eqref{P}, \eqref{defF} and \eqref{stima_ab_uno}, yields that
\begin{align}
	\non
	F'(\ph_\b)\to F'(\ph), \
	\ P(\ph_\b) \to P(\ph), \quad \hbox{both strongly in $C^0(\overline{Q})$.}
\end{align}
Furthermore, we infer that
\Beq
	\non
	B(\ph_{\b}) \to B(\ph)\ \ \hbox{at least strongly in } \L2 H.
\Eeq
It is now a standard matter to complete the details.
Finally, uniqueness follows as a consequence of \cite[Thm.~2.3]{CGRS_ASY}.
\Edim

\dafare{
Let us remark that by \cite[Ex.~2.4]{CGRS_ASY} the authors pointed out
a severe non-uniqueness
result for the system \EQa~if $\a L = 1$.
On the other hand, they pointed out a stability estimate
whenever $\a$ is sufficiently small ($\a < \min \{ \frac 1L, \frac 1{(1+L)^2}\}$), 
which directly implies uniqueness.}
The smallness condition is indeed motivated by the fact that we 
need the uniqueness for system \EQa.
\section{The Control Problem}
\label{SEC_CONTROL_PROBLEM}
\setcounter{equation}{0}

This last section is completely devoted to the investigation of the 
optimal control problem $(CP)$. 
We prove the existence of an optimal control
and point out a variational inequality which characterizes the 
first-order necessary condition for optimality.

\subsection{Existence of Optimal Controls}
To begin with, we check the existence of optimal controls,
namely, we prove Theorem~\ref{THexistenceofcontrol}.

\proof[Proof of Theorem \ref{THexistenceofcontrol}]
For the proof, we employ the direct method of calculus of variations.
In this direction, let us fix a sequence $\graffe{\b_n}_n$ which goes
to zero as $n\to \infty$. Then, let $\graffe {u_n}_n:=\graffe{u_{\b_n}}_n\subset \Uad$ 
be a minimizing sequence for $\J$ which, at every step, is made of
optimal control for $(CP)_{\b_n}$,
and let $(\m_n,\ph_n,\s_n)$ be the corresponding solution to system \EQa.
From the bounds pointed out by estimate \eqref{stima_ab_uno}, we deduce that, as $n \to \infty$,
there exist some $\overline{u} \in \Uad$, a triple $(\bm,\bph,\bs)$, such that,
up to a not relabeled subsequence, the following convergences 
\begin{align}
	\non
	u_n &\to \overline{u} \ \ \hbox{weakly star in } L^\infty(Q)
	\\ \non 
	\m_n &\to \bm  \ \ \hbox{weakly star in } \H1 H \cap \L\infty  V \cap \L2 W   
	\\ \non  
	\ph_n &\to \bph \ \ \hbox{weakly star in } \H1 V \cap \L\infty W 
	\\ \non 
	\s_n &\to \bs \ \ \hbox{weakly star in } \H1 H \cap \L\infty  V \cap \L2 W 
\end{align}
are satisfied.
Moreover, as made in \eqref{strong_phab},
standard compactness results (see, e.g., \cite[Sec. 8, Cor. 4]{Simon})
imply that
\Beq
	\non
	\ph_n \to \bph \ \ \hbox{strongly in } C^0(\overline{Q}),
\Eeq
which also gives sense to the initial condition $\bph(0)=\ph_0$.
Thus, this latter, along with \eqref{P}, \eqref{defF} and \eqref{stima_ab_uno}, allows us
to identify the nonlinear terms in the limit. In fact, as $n \to \infty$, we realize that
\begin{align}
	\non
	F'(\ph_n)\to F'(\bph), \ \ P(\ph_n) \to P(\bph), \quad \hbox{both strongly in $C^0(\overline{Q})$.}
\end{align}
Next, we take into account the variational formulation of 
system \EQa, written for $(\m_n,\ph_n,\s_n)$,
and pass to the limit as $n \to \infty$.
Therefore, we realize that $(\bm,\bph,\bs)$ is the unique solution to \EQa~associated 
with $\overline{u}$. Lastly, invoking the weak sequential lower semicontinuity 
of the cost functional $\J$, it turns out that $\overline{u}$ is the minimizer
we are looking for. 
\qed

\subsection{Approximation of Optimal Controls}
\label{SUB_APPROX}

Once the existence has been proved, we would like 
to characterize the optimality of $(CP)$ on account of some
asymptotic schemes.
If we want to let $\b\searrow 0$ in the optimality condition for $(CP)_\b$
in order to obtain the one for $(CP)$,
we should ensure that every optimal control $\overline u$ 
for $(CP)$ can be found as a limit of a sequence consisting of optimal controls 
for $(CP)_\b$. 
\dafare{
As anticipated, this strong condition is out of reach,
so that we follow a different way making use of the approximated optimal control
problem $(\widetilde{CP})_{\b}$. In fact, we first show that 
$(\widetilde{CP})_{\b}$ can be solved, and then we precise in which sense
it can be useful to deduce the necessary
condition for optimality of $(CP)$. Furthermore,
}
as it complies with the framework of \cite{S}, it is \sfw~to obtain the result below.
\Blem
\label{LM_existencePCtil}
Under the assumptions \tutteleipotesi,
whenever $\b\in(0,1)$ is given, the optimal control problem $(\widetilde{CP})_{\b}$ admits
at least a solution.
\Elem
Moreover, as a consequence of \cite{S},
it also follows the first-order optimality condition for
optimality (compare with Theorem \ref{THMsecondanecab}).
\Bthm
\label{THM_secondanecab_adapted}
Assume that \tutteleipotesi~are satisfied and let $\overline{u}_\b \in \Uad$ 
be an optimal control for $(\widetilde{CP})_{\b}$
with the corresponding optimal state $(\bm_{\b},\bph_{\b},\bs_{\b})$.
Moreover, let $(p_{\b},q_{\b},r_{\b})$ be the solution to the 
adjoint system \EQAggab. Then, the first-order necessary conditions
for optimality reads as follows
\Beq
  \label{secondanecab_adapted}
  \intQ (r_{\b} + \bz \overline{u}_\b + (\overline u_\b - \overline u ))(v - \overline{u}_\b)
  \geq 0 \quad \forall v \in \Uad.
\Eeq
\Ethm

Now, all the ingredients are set and we are in a position to properly 
state the approximation result we mentioned above.
\Bthm
\label{THM_Approximation}
Assume that \tutteleipotesi~are in force. 
Moreover, let $(\bph,\bs,\overline u)$ be an optimal triple for $(CP)$ and
let $\graffe{\b_n}_n$ be a sequence which goes to zero as $n\to \infty$.
Then, there exists an approximating optimal triple, namely
a triple $(\bph_{\b_n}, \bs_{\b_n}, \overline{u}_{\b_n} )$ which solves $(\widetilde{CP})_{\b_n}$
and a not relabeled subsequence such that, as $n\to\infty$, 
we have the following convergences
\begin{align}
	 \label{approx_u}
	\overline u_n:=\overline{u}_{\b_n} &\to \overline u \ \ \hbox{strongly in } L^2(Q)
	\\ 
	 \label{approx_ph}
	\bph_n:=\bph_{\b_n}  &\to \bph \ \ \hbox{weakly star in }  \H1 V \cap \L\infty W
	\\   \label{approx_b}
	\bs_n:= \bs_{\b_n} &\to \bs \ \ \hbox{weakly star in }  \H1 H \cap \L\infty V \cap \L2 W
	\\   \label{approx_contr}
	\widetilde{\J} (\bph_n, \bs_n, \overline{u}_n ) &\to \J(\bph,\bs,\overline u).
\end{align}
\Accorpa\approx approx_u approx_contr
\Ethm
This theorem is the best we can say
as far as the approximation of optimal controls of $(CP)$ by sequences of optimal controls
for an approximating problem is concerned. 
The proof mainly relies on monotonicity and compactness arguments.
\proof[Proof of Theorem \ref{THM_Approximation}]
Let $\b\in(0,1)$, $(\bph_\b, \bs_\b, \overline{u}_\b)$ 
be an optimal triple for $(\widetilde{CP})_\b$, 
which exists by virtue of Lemma \ref{LM_existencePCtil},
and let $\graffe{\b_n}_n$ be a sequence which goes to zero as $n \to\infty$. 
For the sake of simplicity, with $\overline{u}_n$ we denote the optimal 
control associated to 
$\b_n$, that is, $\overline{u}_n:=\overline{u}_{\b_n}$.
Likewise, $\bph_n:=\bph_{\b_n}$ and $\bs_n:=\bs_{\b_n}$.
In view of the boundedness of $\Uad$ and of estimates \accorpa{mab_ma}{sab_sa},
there exist some $\ph, \s, $ and $u$ such that, as $n \to \infty$, the convergences
\Bsist
	&&
	\non
	\overline{u}_n \to {u} \ \ \hbox{weakly star in } L^\infty(Q)
	\\  &&\non
	\bph_n\to \ph \ \ \hbox{weakly star in } {\H1 V \cap \L\infty W}
	\\  &&\non
	\bs_n \to \s \ \ \hbox{weakly star in } \H1 H \cap \L\infty V \cap \L2 W	
\Esist
are verified. Moreover, we also realize that the limit $(\ph,\s,u)$ is an admissible 
triple for $(CP)$.  
Furthermore, we claim that $(\ph,\s,u)$ is nothing but $(\bph,\bs,\overline{u})$,
where $\overline u$ is an optimal control for $(CP)$, whereas $\bph$ 
and $\bs$ are the corresponding states. 
Note that this would imply that the sequence $(\bph_n, \bs_n, \overline u_n)$
approximates $(\bph,\bs,\overline{u})$ in the sense described above.
The weak sequential lower semicontinuity of the adapted cost functional $\widetilde{\J}$ yields that
\begin{align}
	\non &
	\liminf_{n\to\infty} \widetilde{\J} (\bph_n, \bs_n, \overline{u}_n )
	\geq \widetilde{\J} (\ph,\s,u) 
	= \J (\ph,\s,u) +\frac 12 \norma{u-\overline{u}}^2_{L^2(Q)}
	\\ & \quad	
	\label{liminf}
	\geq \J (\bph,\bs,\overline u) +\frac 12 \norma{u-\overline u}^2_{L^2(Q)},
\end{align}
where we also take into account the optimality of
$(\bph,\bs, \overline u)$ for $(CP)$ and the definition of the adpated cost functional
\eqref{adapted_cost}.
On the other hand, the optimality of $(\bph_n, \bs_n, \overline u_n )$ for 
$(\widetilde{CP})_{\b_n}$, implies that
\Beq
	\non
	\widetilde{\J} (\bph_n, \bs_n, \overline u_n )
	\leq \widetilde{\J} (\bph,\bs,\overline u)
	\quad \hbox{ for every $n\in\enne$.}
\Eeq
Hence, passing to the superior limit to both sides leads to deduce that
\Beq
	\label{limsup}
	\limsup_{n\to\infty} \widetilde{\J} (\bph_n, \bs_n, \overline u_n )
	\leq
	\widetilde{\J} (\bph,\bs,\overline u)
	= \J (\bph,\bs,\overline u).
\Eeq
Finally, by combining \eqref{liminf} with \eqref{limsup}, we infer that
\Beq
	\label{andrea3}
	\frac 12 \norma{u- \overline u}^2_{L^2(Q)} = 0.
\Eeq
It is now \sfw~to realize that also the corresponding states coincide,
leading to conclude that $(\ph,\s,u)=(\bph,\bs,\overline{u})$, as we claimed.
Lastly, upon collecting the above information, we have the following chain of
equality
\begin{align}
	\non 
	\J(\bph,\bs,\overline u) 
	&= \widetilde{\J} (\bph,\bs,\overline u)
	= 	\liminf_{n\to\infty} \widetilde{\J} (\bph_n, \bs_n, \overline{u}_n )
	= \limsup_{n\to\infty} \widetilde{\J} (\bph_n, \bs_n, \overline u_n )
	\\ & \label{andrea6}
	= \lim_{n\to\infty} \widetilde{\J} (\bph_n, \bs_n, \overline u_n )
	= \lim_{n\to\infty} {\J} (\bph_n, \bs_n, \overline u_n )
	+ \frac 12 \norma{\overline u_n - \uopt}^2_{L^2(Q)},
\end{align}
Thus, we are reduced to show \eqref{approx_u}. Up to now, we have just proven that
the weak limit of ${\overline u}_n$ is $\uopt$. On the other hand, 
it easily follows from the above estimates, along with the lower semicontinuity of the cost functional,
that 
\begin{align*}
		\non
		\J(\bph,\bs,\overline u) 
		& \leq \liminf_{n\to\infty} \J(\bph_n, \bs_n, \overline u_n) 
		\leq \limsup_{n\to\infty} \J(\bph_n, \bs_n, \overline u_n) 
		\\ & 
		\leq \limsup_{n\to\infty} \widetilde{\J} (\bph_n, \bs_n, \overline u_n) 
		= \lim \widetilde{\J} (\bph_n, \bs_n, \overline u_n) 
		= \J(\bph,\bs,\overline u),
\end{align*}
so that
\begin{align}
	\label{andrea5}
	\J(\bph,\bs,\overline u)
	=
	\lim_{n\to\infty} \J(\bph_n, \bs_n, \overline u_n)
\end{align}
is verified. Therefore, by combining \eqref{andrea6} with \accorpa{andrea3}{andrea5},
we deduce that
\begin{align*}
	\frac 12 \norma{\overline{u}_n  - \uopt}^2_{L^2(Q)} \to 0,
\end{align*}
which conclude the proof.
\qed
\subsection{The Adjoint System}
\label{ADJOINT}
\setcounter{equation}{0}

Here, we are going to investigate the adjoint system
proving Theorem \ref{THM_stimeaggiuntoab}.
In order to avoid a heavy notation, we will omit 
writing the subscript $\b$ on the variables which occur in the calculations below,
while we will reintroduce the correct notation at the end of each estimate.
Before moving on, let us set
\Beq
	\non
	Q_t^T:= \Omega \times [t, T], \quad \hbox{for every $t \in [0,T).$}
\Eeq
Below, we will proceed quite formally. The justification can be carried out rigorously, 
e.g., introducing a Galerkin scheme.
Moreover, the adjoint system is linear and therefore, the
uniqueness part easily follows by applying standard arguments from
the existence part. On the other hand, system \EQAggab~has
already been studied in \cite[Sec.~4.4]{S}, and we refer there the interested
reader for the details of the Galerkin technique.
\proof[Proof of Theorem \ref{THM_stimeaggiuntoab}]
The a priori estimates we are going to point out
will allow us to justify in a rigorous way the passage $\b\searrow 0$ in system \EQAggab.

{\bf First estimate}
To begin with, we add to both sides of \eqref{EQAggabgseconda} the term $p$. Then, we test
\eqref{EQAggabprima} by $-q$, this new second equation by $-\dt p$, and \eqref{EQAggabterza} by $r$.
Summing up and integrating over $Q_t^T$ lead to
\begin{align}
	\non
	&
	\frac \b2 \IO2 q
	+ \intQtT \dt p \, q
	+ \intQtT |{\nabla q}|^2
	+ \intQtT F''(\bph) |{q}|^2
	+ \frac 12 \IO2 {p}
	\\ \non & \qquad
	- \intQtT \dt p \, q
	+\a \intQtT |\dt p|^2
	+ \frac 12 \IO2 {\nabla p}
	+ \frac 12 \IO2 {r}
	+ \intQtT |{\nabla r}|^2
	\\ \non & \quad
	=
	 \frac \b2 \iO |{q(T)}|^2
	+ \frac 12 \iO | {p(T)} |^2
	+ \frac 12 \iO |{\nabla p(T)}|^2
	+ \frac 12 \iO | {r(T)} |^2
	\\ \non & \qquad
	- \bQ \intQtT (\bph-\phQ)q
	+ \bQh \intQtT (\bs-\sQ)r
	+ \intQtT P'(\bph)(\bs-\bm)(r-p)q
	\\ & \qquad
	+ \intQtT P(\bph)(p-r) \dt p
	- \intQtT p \, \dt p	
	- \intQtT P(\bph)(r-p)r.
	\label{andrea1}
\end{align}
On the \lhs~two integrals cancel out and, despite the fourth term,
the others are nonnegative. As the fourth term is concerned, we remind 
that $F''$ is bounded below in terms
of the \Lip\ constant $L$, so that we have 
\Beq
	\non
	\intQtT F''(\bph) |{q}|^2 \geq - L \intQtT |q|^2.
\Eeq
Next, we test \eqref{EQAggabgseconda} by $K q$,
for a positive constant $K$ yet to be determined, 
and integrate over $Q_t^T$ to get
\Bsist
	K \intQtT |q|^2
	= \a K \intQtT \dt p \, q 
	- K \intQtT \nabla p \cdot \nabla q
	- K \intQtT P(\bph)(p-r)q.
	\label{andrea2}
\Esist
Then, after making use of the definition of the final conditions \eqref{ICEQAggab},
we add \accorpa{andrea1}{andrea2} to obtain that
\begin{align}
	\non &
	\frac \b2 \IO2 q
	+ (K-L) \intQtT |q|^2
	+ \intQtT |{\nabla q}|^2
	+\a \intQtT |\dt p|^2
	\\ \non & \qquad
	+ \frac 12 \IO2 {p}		
	+ \frac 12 \IO2 {\nabla p}
	+ \frac 12 \IO2 {r}
	+ \intQtT |{\nabla r}|^2
	\\ \non & \quad
	\leq
	\frac {\bOh}2 \iO |\bs(T)-\sO|^2
	- \bQ \intQtT (\bph-\phQ)q
	+ \bQh \intQtT (\bs-\sQ)r
	\\ \non & \qquad
	+ \intQtT P'(\bph)(\bs-\bm)(r-p)q
	+ \intQtT P(\bph)(p-r) \dt p
	\\ \non & \qquad
	- \intQtT p \, \dt p	
	- \intQtT P(\bph)(r-p)r
	+ \a K \intQtT \dt p \, q 
	\\ \non & \qquad
	- K \intQtT \nabla p \cdot \nabla q
	- K \intQtT P(\bph)(p-r)q,
\end{align}
where we denote by $I_1,\dots,I_{10}$ the integrals on the \rhs, in that order.
Now, we start estimating the terms on the \rhs.
Owing to assumptions \accorpa{constants}{targets}, we easily realize that
\Bsist
	\non
	&&  |I_1|
	\leq
	c.
\Esist
Meanwhile, the integrals $I_2$ and $I_3$ can be easily managed by applying the Young inequality
and the fact that $\bph$ and $\bs$, as solutions to \EQab, satisfy \eqref{stima_ab_uno}.
In fact, we have that
\Bsist
	\non
	&& |I_2|  + |I_{3}|   
	\leq 
	\d \intQtT |q|^2
	+ \frac 12 \intQtT |r|^2 
	+ \cd,
\Esist
for a small and positive $\d$ yet to be determined.
Invoking the \Holder~and Young inequalities, 
the continuous embeddings $V\subset\Lx4$ and $V\subset\Lx6$, assumption \eqref{P},
and estimate \eqref{stima_ab_uno},  
we infer that
\begin{align}
	\non
 	|I_4|  
	&\leq
	c \intQtT |\bs - \bm | | r-p | |q |
	\leq c \inttt \norma{\bs - \bm }_{6} \norma{r-p}_{2}\norma{q}_{3}
	\\ \non & 
	\leq \d \inttt \normaV{q}^2 + \cd \inttt (\normaV{\bs}^2 
	+ \normaV{\bm}^2) (\normaH{r}^2+\normaH{p}^2)
	\\ \non & 
	\leq \d \intQtT (|q|^2+|\nabla q|^2) + \cd \intQtT (|r|^2+|p|^2).
\end{align}
By the same token, using \eqref{young}, we get
\Bsist
	\non
	&& |I_5|  
	\leq
	c \intQtT | p-r| | \dt p |
	\leq 
	\d \intQtT |\dt p|^2 + \cd \intQtT (|p|^2+|r|^2).
\Esist
Using the Young inequality once more, we realize that
\Bsist
	\non
	&& |I_{6}| +|I_{7}|  
	\leq
	\d \intQtT |\dt p|^2 + \cd \intQtT |p|^2
	+ c \intQtT |r|^2,
\Esist
and also that
\Bsist
	\non
	&& |I_{9}|  + |I_{10}| 
	\leq
	\d \intQtT |\nabla q|^2 + \cd \intQtT |\nabla p|^2
	+ \d \intQtT |q|^2 + \cd \intQtT (|p|^2+|r|^2).
\Esist
Lastly, owing to \eqref{young}, $I_{8}$ can be dealt by
\Bsist
	\non
	&& |I_{8}|  
	= \Big| \a K \intQtT \dt p \, q \Big|
	\leq 
	\frac K2 \intQtT | q |^2
	+ \frac {\a^2 K}{2} \intQtT | \dt p |^2.
\Esist
Collecting all the previous estimates, we realize that the backward-in-time
Gronwall lemma yields the estimate we are looking for, provided we check that
$K$ and $\d$ can be chosen in such a way to satisfy the following
condition
\Beq
	\non
	\min \Bigl\{ K-\frac K2 - L -3 \d, 1 - 2\d, \a - \frac{\a^2 K}{2} - 2\d
	\Bigr \}>0.
\Eeq
Actually, considering that $\d$ can be taken arbitrarily small, we are reduced to show that
there exists a positive constant $K$ such that
\Beq
	\non
	\min \Bigl\{\frac K2 - L , \a - \frac{\a^2 K}{2} \Bigr \}>0.
\Eeq
\dafare{
Let us claim that this is satisfied if $\a$ is small enough.
For instance, we take $K=3L$, so that 
\Beq
	\non
	\min \Bigl\{\frac K2 - L , \a - \frac{\a^2 K}{2} \Bigr \}=
	\min \Bigl\{\frac L2 , \frac {\a}2 (2-\a3L) \Bigr \}.
\Eeq
Hence, we assume at once $\a$ to be small in order that $2-\a3L>0$.} 
Lastly, we pick $\d$ 
small enough to conclude.
Finally, a Gronwall argument implies that
\begin{align}
	& \non
	\b^{1/2} \norma {\qab}_{\L\infty H}
	+ \norma{\qab}_{\L2 V}
	+ \norma{\pab}_{\H1 H \cap \L\infty V}
	\\   & \quad
	+ \norma{\rab}_{\L\infty H \cap \L2 V}
	\leq c.
\end{align}

{\bf Second estimate}   
Then, we test \eqref{EQAggabgseconda}
by $- \Delta p$ and integrate over $Q_t^T$ to obtain that
\Bsist
	\non &&
	\frac \a2 \iO |\nabla p(t)|^2
	+ \intQtT |\Delta p|^2
	= \frac \a2 \iO |\nabla p(T)|^2
	+ \intQtT q \, \Delta p
	+ \intQtT P(\bph)(p-r)\Delta p,
\Esist
where we denote the terms on the \rhs~by $I_1,I_2,$ and $I_3$, respectively.
Owing to the final conditions \eqref{ICEQAggab}, we easily conclude that $I_1=0$.
Moreover, Young's inequality, combined with the boundedness of $P$, gives that
\Bsist
	&&
	\non
	|I_2|+|I_3|
	\leq 
	\frac 12 \intQtT |\Delta p|^2
	+ \intQtT |q|^2
	+ c \intQtT (|p|^2 +|r|^2 ).
\Esist
Hence, the previous estimate produces 
\Beq
	\norma{\nabla \pab}_{\L\infty H}
	+ \norma{\Delta \pab}_{\L2 H}
	\leq
	c,
\Eeq
from which, applying standard elliptic regularity results for homogeneous Neumann
boundary problems, we obtain that
\Beq
	\norma{\pab}_{\L\infty V \cap \L2 W}
	\leq
	c.
\Eeq

{\bf Third estimate}
Furthermore, we test \eqref{EQAggabterza} first by $-\dt r$ and 
secondly by $-\Delta r$ to get the following parabolic regularity
\Beq
	\norma{\rab}_{\H1 H \cap \L\infty V \cap \L2 W} \leq c.
\Eeq

{\bf Fourth estimate}  
Next, by testing \eqref{EQAggabprima} by $\Delta q$ and integrating over
$Q_t^T$, we find that
\Bsist
	\non &&
	\frac \b2 \iO |\nabla q(t)|^2
	+ \intQtT |\Delta q|^2
	= 	\frac \b2 \iO |\nabla q(T)|^2
	+ \intQtT \dt p \, \Delta q
	+ \intQtT F''(\bph)q \, \Delta q
	\\ \non && \quad
	- \intQtT P'(\bph)(\bs-\bm)(r-p)\Delta q
	+ \intQtT \bQ (\bph-\phQ)\Delta q,
\Esist
where we indicate the terms on the \rhs~by $I_1, \dots, I_5$, in that order.
In a similar fashion as in the previous estimates, we first observe that 
$I_1=0$. Next, owing to the
Young and \Holder\ inequalities, to the boundary of $P$, 
and to the continuous embedding \eqref{VinL6},
the remaining integrals can be dealt as
\begin{align}
 	\non
 	& |I_2| + |I_3| + |I_4| + |I_5|
 	\leq
 	\frac 45 \intQtT |\Delta q|^2
	+ c \intQtT |\dt p|^2  	
 	+ c \intQtT |q|^2
 	+c \intQtT (|p|^2 +|r|^2 ) + c,
\end{align}
where estimate \eqref{stima_ab_uno} for the
solutions $\bm$ and $\bs$, is also taken into account. 
Hence, we have that 
\Beq
	\b^{1/2}\norma{\nabla \qab}_{\L\infty H}
	+ \norma{\Delta \qab}_{\L2 H}
	\leq
	c,
\Eeq
and the regularity results for elliptic equations with 
homogeneous Neumann boundary conditions, entails that
\Beq
	\b^{1/2}\norma{\nabla \qab}_{\L\infty H}
	+ \norma{\qab}_{\L2 W}
	\leq
	c.
\Eeq
%
{\bf Fifth estimate}
Lastly, we rearrange equation \eqref{EQAggabprima} in the following way
\Beq
	\non
	\b \dt q
	=
	\dt p
	- \Delta q 
	+ F''(\bph)q
  	- P'(\bph)(\bs - \bm)(r-p)
  	+ \bQ(\bph - \phQ).
\Eeq
Therefore, by comparison in the above equation, we also realize that
\Beq
	\label{norma_star}
	\b \norma{\dt q_\b}_{\L2 H} \leq c.
\Eeq

{\bf Passing to the limit}
Summing up, upon combining the above estimates we recover estimate 
\eqref{stimaagg}. Moreover, we infer that
there exist some variables $q,p$ and $r$ such that,
up to a not relabeled subsequence, as $\b \searrow 0$, the convergences
mentioned by \accorpa{q_abtoqa}{bqtozero_due} hold.
Furthermore, these uniform bounds, along with 
standard compactness embedding results,
allow us to recover also the following strong convergences
\Bsist
	&&	
	\label{pstrong}
	\pab \to p \quad \hbox{strongly in } \C0 H \cap \L2 V
	\\  && 
	\label{rstrong}
	\rab \to r \quad \hbox{strongly in } \C0 H \cap \L2 V.
\Esist
Then, we try to draw some consequences from the aforementioned a priori bounds 
in order to pass to the limit, as $\b\searrow 0$, in the adjoint system \EQAggab. 
For the sake of convenience, we 
rewrite its variational formulation which can be 
obtained by testing the system by
an arbitrary $v \in V$ and integrating over $\Omega$. It reads as follows
\begin{align}
	\non &
	 \b \iO\dt q_{\b}(t) \, v	
	- \iO  \dt p_{\b}(t) \, v
	- \iO \nabla q_{\b}(t) \cdot \nabla v
	- \iO  F''(\bph_{\b}(t))q_{\b}(t) \, v
	\\ & \qquad \non
	+ \iO P'(\bph_{\b}(t))(\bs_{\b}(t)-\bm_{\b}(t))(r_{\b}(t)-p_{\b}(t)) v
	= \iO \bQ (\bph_{\b}(t)-\phQ(t)) v
	\\ \non &
	\hspace{7cm}
	\quad \hbox{for every $v\in V$, \  $\aat$}
	\\ & \non
	\iO q_{\b}(t) v
	{-\a}\iO \dt p_{\b}(t) \, v
	{+} \iO \nabla p_{\b}(t) \cdot \nabla v
	\\ \non &
	\hspace{2cm}
	+ \iO P(\bph_{\b}(t))(p_{\b}(t) -r_{\b}(t)) v = 0
	\quad \hbox{for every $v\in V$, \  $\aat$}
	\\ 	& \non
	- \iO \dt r_{\b}(t) \,v
	+ \iO \nabla r _{\b}(t)\cdot \nabla v
	+ \iO P(\bph_{\b}(t))(r_{\b}(t)-p_{\b}(t)) v
	\\ \non &
	\hspace{4cm}
	= \iO \bQh (\bs_{\b}(t)-\sQ(t)) v
		\quad \hbox{for every $v\in V$, \  $\aat$,}
\end{align}
and the final condition
\Beq
	\non
	\iO r_\b(T) v = \iO \bOh (\bs_\b (T)- \sO) v 
	\quad
	\hbox{for every $v\in V$.}
\Eeq
At this point, we would invoke the above convergences \convadjab~to 
show that in the limit, as $\b \searrow 0$, we find
\begin{align}
	\non &
	- \iO  \dt p(t) \, v
	- \iO \nabla q(t) \cdot \nabla v
	- \iO  F''(\bph(t))q(t) \, v
	\\ & \quad \non
	+ \iO P'(\bph(t))(\bs(t)-\bm(t))(r(t)-p(t)) v
	= \iO \bQ (\bph(t)-\phQ(t)) v
	\\ \non &
	\hspace{7cm}
	\quad \hbox{for every $v\in V$, \  $\aat$}
	\\ & \non
	\iO q(t) v
	{-\a}\iO \dt p(t) \, v
	{+} \iO \nabla p(t) \cdot \nabla v
	\\ \non &
	\hspace{2.5cm}
	+ \iO P(\bph(t))(p(t) -r(t)) v 	
	= 0 \quad \hbox{for every $v\in V$, \  $\aat$}
	\\ 	& \non
	- \iO \dt r(t) \,v
	+ \iO \nabla r (t)\cdot \nabla v
	+ \iO P(\bph_{\b}(t))(r_{\b}(t)-p_{\b}(t)) v
	\\ \non &
	\hspace{4cm}
	= \iO \bQh (\bs(t)-\sQ(t)) v
	\quad \hbox{for every $v\in V$, \  $\aat$,}
\end{align}
and the final condition
\Beq
	\non
	\iO r(T) v = \iO \bOh (\bs (T)- \sO) v 
	\quad
	\hbox{for every $v\in V$,}
\Eeq
which corresponds to the variational formulation associated to system \EQAgga.
Nonetheless, since there appear some nonlinear terms, some care is in order. 
First, let us recall that both $P$ and $F$ are regular due to \eqref{P} and \eqref{regpot}
and that \eqref{VIIest_b} holds.
Hence, exploiting the strong convergence \eqref{strong_phab}, 
we claim that, as $\b\searrow 0$, we have
\begin{align}
	 \label{stronF}
	F''(\bph_{\b}) &\to F''(\bph) \quad  \hbox{strongly in } C^0(\overline{Q})
	\\ 
	\label{stronP}
	P(\bph_{\b}) &\to P(\bph)
	\quad \hbox{strongly in } C^0(\overline{Q}).
\end{align}
\Accorpa\strongFP  {stronF} {stronP} 
To prove the former, it suffices to combine \eqref{growth} and \eqref{strong_phab} with 
the estimate \eqref{VIIest_b}, while for
the latter we simply account for \eqref{strong_phab} and for 
the boundedness of $P$.
 Moreover, having in mind the weak convergence \eqref{q_abtoqa} and 
the strong ones \accorpa{pstrong}{rstrong}, we can prove that
the nonlinear terms can be identified in the limit.
In fact, from \eqref{q_abtoqa}, we have that 
\Beq
	\non
	\qab \to q \quad  \hbox{at least weakly in } \L2 H,
\Eeq
which, along with \eqref{stronF}, leads to infer that
\Beq
	\non
	F''(\bph_{\b})q_{\b} \to F''(\bph)q
	\quad 
	\hbox{at least weakly in $\L2 H$.}
\Eeq
Similarly, combining \accorpa{pstrong}{rstrong} with \eqref{stronP}, we also deduce that
\begin{align}
	\non
	P(\bph_{\b})(p_{\b}-r_{\b}) &\to P(\bph)(p -r)
	\quad 
	\hbox{strongly in $L^2(Q)$}
	\\ 
	\non
	P'(\bph_{\b})(\bs_{\b}-\bm_{\b})(r_{\b}-p_{\b}) 
	&\to P'(\bph)(\bs-\bm)(r-p)
	\quad 
	\hbox{strongly in $L^1(Q)$},
\end{align}
where the first one requires the help of the strong convergence
pointed out by \eqref{strong_phab} and the boundedness of $P$, whereas
in the second we
again owe to \eqref{strong_phab}, along with the strong convergences 
\accorpa{strong_mab}{strong_sab}. 
To completely recover system \EQAgga~it suffices to check
that the regularity is enough to rewrite the system in a strong form.
So, this is the sense in which we can say that system \EQAggab~converges, as $\b\searrow 0$, to 
\EQAgga.
\qed

\subsection{First-order Necessary Condition}
We conclude the paper providing the first-order necessary condition
that an optimal control, which exists in view of Theorem \ref{THexistenceofcontrol},
has to verify. 

\proof[Proof of Theorem \ref{THMsecondanec}]
As previously mentioned, in order to get inequality \eqref{secondanec}, we
cannot pass to the limit as $\b\searrow 0$ in the variational inequality 
\eqref{secondanecab} since nothing ensures that in such a passage,
the control $\overline{u}_\b$ will converge to a limit that is also optimal for $(CP)$.
Therefore, the investigation made in the subsection \ref{SUB_APPROX} helps to rigorously
handle this issue. Indeed, we are going to consider a sequence $\graffe{\b_n}$
which goes to zero as $n\to \infty$ and then we take into account 
$\overline{u}_n:= \overline{u}_{\b_n}$ instead of $\overline{u}_\b$. 
Therefore, after extraction of a subsequence $\graffe{\b_{n_k}}$, 
the asymptotics pointed out by \convadjab~and \approx~allow 
us to pass to the limit as $k\to \infty$
in \eqref{secondanecab_adapted} to obtain \eqref{secondanec}.

Furthermore, the last sentence immediately follows as a \sfw~application
of the Hilbert projection theorem, since $\Uad$ is a non-empty, closed and convex subset 
of $\L2 H$.
Moreover, let us note that \eqref{secondanec} implies that, 
whenever $\bz > 0$, the optimal control $\overline{u}$
can be implicitly characterized as follows (see, e.g., \cite{Trol})
\Beq
	\non
	\overline{u}(x,t)=\max \bigl\{ 
	u_*(x,t), \min\graffe{u^*(x,t),-\frac 1{\bz} r(x,t)} 
	\bigr\} \quad \aaQ .
\Eeq
\qed


\subsection*{Acknowledgments}
The author wishes to thanks Professor Pierluigi Colli for several useful 
discussions and suggestions which have improved the manuscript.
%
%
%
\vspace{3truemm}

\Begin{thebibliography}{10}
\footnotesize

\bibitem{Agosti}
A. Agosti, P.F. Antonietti, P. Ciarletta, M. Grasselli and M. Verani,
A Cahn-Hilliard-type equation with application to tumor growth dynamics, 
{\it Math. Methods Appl. Sci.}, {\bf 40} (2017), 7598-–7626.

\bibitem{BARBU}
V. Barbu,
Necessary conditions for nonconvex distributed control problems governed by 
elliptic variational inequalities, 
{\it J. Math. Anal. Appl.} {\bf 80} (1981), 566-597.

\bibitem{BAR}
V. Barbu, 
``Nonlinear semigroups and differential equations in Banach spaces",
Noordhoff International Publishing, Leyden, 1976.

\bibitem{BRZ}
H. Brezis,
``Op\'erateurs maximaux monotones et semi-groupes de contractions dans les
espaces de Hilbert'', North-Holland Math. Stud. {\bf 5}, North-Holland, Amsterdam, 1973.

\bibitem{CRW}
C. Cavaterra, E. Rocca and H. Wu,
Long-time Dynamics and Optimal Control of a Diffuse Interface Model for Tumor Growth,
{\it preprint arXiv:1901.07500 [math.AP],} (2018), 1-–36.

\bibitem{CFGS}
P. Colli, M.H. Farshbaf-Shaker, G. Gilardi and J. Sprekels,
Optimal boundary control of a viscous Cahn–-Hilliard system
with dynamic boundary condition and double obstacle potentials,
{\it SIAM J. Control Optim.} {\bf 53} (2015), 2696-2721.
\bibitem{CFS}
P. Colli, M.H. Farshbaf-Shaker and J. Sprekels,
A deep quench approach to the optimal control of an Allen–-Cahn equation
with dynamic boundary conditions and double obstacles,
{\it Appl. Math. Optim.} {\bf 71} (2015), 1-24.

\bibitem{CGH}
P. Colli, G. Gilardi and D. Hilhorst,
On a Cahn-Hilliard type phase field system related to tumor growth,
{\it Discrete Contin. Dyn. Syst.} {\bf 35} (2015), 2423-2442.

\bibitem{CGMR_sing}
P. Colli, G. Gilardi, G. Marinoschi and E. Rocca,
Optimal control for a phase field system with a possibly singular potential,
{\it Math. Control Relat. Fields\/} {\bf 6} (2016), 95-112.

\bibitem{CGMR_cons}
P. Colli, G. Gilardi, G. Marinoschi and E. Rocca,
Optimal control for a conserved phase field system with a possibly singular potential,
{\it Evol. Equ. Control Theory\/} {\bf 7} (2018), 95-116.

\bibitem{CGRS_VAN}
P. Colli, G. Gilardi, E. Rocca and J. Sprekels,
Vanishing viscosities and error estimate for a Cahn–Hilliard type phase field system related to tumor growth,
{\it Nonlinear Anal. Real World Appl.} {\bf 26} (2015), 93-108.

\bibitem{CGRS_OPT}
P. Colli, G. Gilardi, E. Rocca and J. Sprekels,
Optimal distributed control of a diffuse interface model of tumor growth,
{\it Nonlinearity} {\bf 30} (2017), 2518-2546.

\bibitem{CGRS_ASY}
P. Colli, G. Gilardi, E. Rocca and J. Sprekels,
Asymptotic analyses and error estimates for a \CH~type phase field system modeling tumor growth,
{\it Discrete Contin. Dyn. Syst. Ser. S} {\bf 10} (2017), 37-54.

\bibitem{CGS}
P. Colli, G. Gilardi and J. Sprekels,
On the \CH~equation with dynamic 
boundary conditions and a dominating boundary potential,
{\it J. Math. Anal. Appl.\/} {\bf 419} (2014) 972-994.

\bibitem{CGS_OPT}
P. Colli, G. Gilardi and J. Sprekels,
A boundary control problem for the viscous \CH~equation 
with dynamic boundary conditions,
{\it Appl. Math. Opt.\/} {\bf 73} (2016) 195-225.


\bibitem{CGS_nonst}
P. Colli, G. Gilardi and J. Sprekels,
Optimal boundary control of a nonstandard viscous \CH~system with dynamic boundary condition,
{\it Nonlinear Anal.\/} {\bf 170} (2018), 171-196.

\bibitem{CGS14}
P. Colli, G. Gilardi and J. Sprekels,
Optimal velocity control of a viscous Cahn-Hilliard system with convection and dynamic boundary conditions, 
{\it SIAM J. Control Optim.\/} {\bf 56} (2018), 1665-1691.

\bibitem{CGS_DQ}
P. Colli, G. Gilardi and J. Sprekels,
Optimal velocity control of a convective Cahn-–Hilliard system with double obstacles
and dynamic boundary conditions: a `deep quench' approach.
{\it J. Convex Anal.}, to appear (2018).
%

\bibitem{CS}
P. Colli and J. Sprekels,
Optimal control of an Allen-Cahn equation 
with singular potentials and dynamic boundary condition,
{\it SIAM J. Control Optim.\/} {\bf 53} (2015) 213-234.


\bibitem{CLLW}
V. Cristini, X. Li, J.S. Lowengrub, S.M. Wise,
Nonlinear simulations of solid tumor growth using a mixture model: invasion and branching.
{\it J. Math. Biol.} {\bf 58} (2009), 723–-763.

\bibitem{CL}
V. Cristini, J. Lowengrub,
Multiscale Modeling of Cancer: An Integrated Experimental and Mathematical
{\it Modeling Approach. Cambridge University Press}, Leiden (2010).

\bibitem{DFRGM}
M. Dai, E. Feireisl, E. Rocca, G. Schimperna, M. Schonbek,
Analysis of a diffuse interface model of multispecies tumor growth,
{\it Nonlinearity\/} {\bf  30} (2017), 1639.

\bibitem{EK_ADV}
M. Ebenbeck and P. Knopf,
Optimal control theory and advanced optimality conditions for a diffuse interface model of tumor growth 
{\it preprint arXiv:1903.00333 [math.OC],} (2019), 1-34.

\bibitem{EK}
M. Ebenbeck and P. Knopf,
Optimal medication for tumors modeled by a Cahn-Hilliard-Brinkman equation,
{\it preprint arXiv:1811.07783 [math.AP],} (2018), 1-26.

\bibitem{EGAR}
M. Ebenbeck and H. Garcke,
Analysis of a Cahn–-Hilliard–-Brinkman model for tumour growth with chemotaxis.
{\it J. Differential Equations,} (2018) https://doi.org/10.1016/j.jde.2018.10.045.

\bibitem{FGR}
S. Frigeri, M. Grasselli, E. Rocca,
On a diffuse interface model of tumor growth,
{\it  European J. Appl. Math.\/} {\bf 26 } (2015), 215-243. 

\bibitem{FLRS}
S. Frigeri, K.F. Lam, E. Rocca, G. Schimperna,
On a multi-species Cahn-Hilliard-Darcy tumor growth model with singular potentials,
{\it Comm. in Math. Sci.} {\bf  (16)(3)} (2018), 821-856. 

\bibitem{FRL}
S. Frigeri, K.F. Lam and E. Rocca,
On a diffuse interface model for tumour growth with non-local interactions and degenerate 
mobilities,
In {\sl  Solvability, Regularity, and Optimal Control of Boundary Value Problems for PDEs},
P. Colli, A. Favini, E. Rocca, G. Schimperna, J. Sprekels (ed.),
{\it Springer INdAM Series,} {\bf 22}, Springer, Cham, 2017.
%
\bibitem{GARL_3}
H. Garcke and K. F. Lam,
Global weak solutions and asymptotic limits of a Cahn-Hilliard-Darcy system modelling tumour growth,
{\it AIMS Mathematics} {\bf 1 (3)} (2016), 318-360.
\bibitem{GARL_1}
H. Garcke and K. F. Lam,
Well-posedness of a Cahn-–Hilliard–-Darcy system modelling tumour
growth with chemotaxis and active transport,
{\it European. J. Appl. Math.} {\bf 28 (2)} (2017), 284-316.
\bibitem{GARL_2}
H. Garcke and K. F. Lam,
Analysis of a Cahn-Hilliard system with non-zero Dirichlet 
conditions modeling tumor growth with chemotaxis,
{\it Discrete Contin. Dyn. Syst.} {\bf 37 (8)} (2017), 4277-4308.
\bibitem{GARL_4}
H. Garcke and K. F. Lam,
On a Cahn-Hilliard-Darcy system for tumour growth with solution dependent source terms, 
in {\sl Trends on Applications of Mathematics to Mechanics}, 
E.~Rocca, U.~Stefanelli, L.~Truskinovski, A.~Visintin~(ed.), 
{\it Springer INdAM Series} {\bf 27}, Springer, Cham, 2018, 243-264.
\bibitem{GAR}
H. Garcke, K. F. Lam, R. N\"urnberg and E. Sitka,
A multiphase Cahn-Hilliard-Darcy model for tumour growth with necrosis,
{\it Mathematical Models and Methods in Applied Sciences} {\bf 28 (3)} (2018), 525-577.

\bibitem{GARLR}
H. Garcke, K. F. Lam and E. Rocca,
Optimal control of treatment time in a diffuse interface model of tumor growth,
{\it Appl. Math. Optim.} {\bf 78}(3) (2018), {495-544}.
\bibitem{GLSS}
H. Garcke, K.F. Lam, E. Sitka, V. Styles,
A Cahn–-Hilliard-–Darcy model for tumour growth with chemotaxis and active transport.
{\it Math. Models Methods Appl. Sci. } {\bf 26(6)} (2016), 1095-–1148.

\bibitem{GiMiSchi} 
G. Gilardi, A. Miranville and G. Schimperna,
On the \CH~equation with irregular potentials and dynamic boundary conditions,
{\it Commun. Pure Appl. Anal.\/} {\bf 8} (2009) 881-912.

\bibitem{OHP}
A. Hawkins, J.T Oden, S. Prudhomme,
General diffuse-interface theories and an approach to predictive tumor growth modeling. 
{\it Math. Models Methods Appl. Sci.} {\bf 58} (2010), 723–-763. 

\bibitem{HDPZO}
A. Hawkins-Daarud, S. Prudhomme, K.G. van der Zee, J.T. Oden,
Bayesian calibration, validation, and uncertainty quantification of diffuse 
interface models of tumor growth. 
{\it J. Math. Biol.} {\bf 67} (2013), 1457–-1485. 

\bibitem{HDZO}
A. Hawkins-Daruud, K. G. van der Zee and J. T. Oden, Numerical simulation of
a thermodynamically consistent four-species tumor growth model, Int. J. Numer.
{\it Math. Biomed. Engng.} {\bf 28} (2011), 3–-24.

\bibitem{HKNZ}
D. Hilhorst, J. Kampmann, T. N. Nguyen and K. G. van der Zee, Formal asymptotic
limit of a diffuse-interface tumor-growth model, 
{\it Math. Models Methods Appl. Sci.} {\bf 25} (2015), 1011-–1043.

\bibitem{Kur}
S. Kurima,
Asymptotic analysis for Cahn-Hilliard type phase field systems related to tumor 
growth in general domains,
{\it Math. Methods in the Appl. Sci} {} (2019), 
https://doi.org/10.1002/mma.5520.

\bibitem{LDY}
O.A. Lady\v zenskaja, V.A. Solonnikov and N.N. Uralceva,
``Linear and quasilinear equations of parabolic type'',
Mathematical Monographs Volume {\bf 23},
{\it American mathematical society}, Providence, 1968.

\bibitem{Lions}
J.-L. Lions,
``\'Equations diff\'erentielles op\'erationnelles et probl\`emes aux limites'',
Grundlehren, Band~111,
Springer-Verlag, Berlin, 1961.


\bibitem{Mir_CH}
A. Miranville,
The Cahn-Hilliard equation and some of its variants, 
{\it AIMS Mathematics,} {\bf 2} (2017), 479–-544.

\bibitem{MRS}
A. Miranville, E. Rocca, and G. Schimperna,
On the long time behavior of a tumor growth model,
{\it J. Differential Equations,} (2019), {https://doi.org/10.1016/j.jde.2019.03.028.}.

\bibitem{MirZelik}
A. Miranville and S. Zelik,
Robust exponential attractors for \CH~type equations with singular potentials,
{\it Math. Methods Appl. Sci.\/} {\bf 27} (2004) 545-582.

\bibitem{S_a}
A. Signori,
Vanishing parameter for an optimal control problem modeling tumor growth.
{\it Preprint: arXiv:1903.04930 [math.AP]} (2019), 1-22.

\bibitem{S_DQ}
A. Signori,
Optimality conditions for an extended tumor growth model with 
double obstacle potential via deep quench approach. 
{\it Preprint: arXiv:1811.08626 [math.AP]} (2018), 1-25.

\bibitem{S}
A. Signori,
Optimal distributed control of an extended model of tumor 
growth with logarithmic potential. 
{\it Appl. Math. Optim.} (2018), https://doi.org/10.1007/s00245-018-9538-1.

\bibitem{Simon}
J. Simon,
{Compact sets in the space $L^p(0,T; B)$},
{\it Ann. Mat. Pura Appl.\/} 
{\bf 146~(4)} (1987) 65-96.

\bibitem{Trol}
F. Tr\"oltzsch,
Optimal Control of Partial Differential Equations. Theory, Methods and Applications,
{\it Grad. Stud. in Math.,} Vol. {\bf 112}, AMS, Providence, RI, 2010.

\bibitem{WLFC}
S.M. Wise, J.S. Lowengrub, H.B. Frieboes, V. Cristini,
Three-dimensional multispecies nonlinear tumor growth—I: model and numerical method. 
{\it J. Theor. Biol.} {\bf 253(3)} (2008), 524–-543.

\bibitem{WZZ}
X. Wu, G.J. van Zwieten and K.G. van der Zee, Stabilized second-order splitting
schemes for \CH~models with applications to diffuse-interface tumor-growth models, 
{\it Int. J. Numer. Meth. Biomed. Engng.} {\bf 30} (2014), 180-203.

\End{thebibliography}

\End{document}

\bye